\documentclass[10pt]{article}
\usepackage{amsmath}
\usepackage{latexsym}
\usepackage{amssymb}
\usepackage{amsfonts}
\usepackage{amsthm}
\title{SPENCER OPERATOR AND APPLICATIONS:\\
From Continuum Mechanics to Mathematical Physics}
\author{J.F. Pommaret \\ CERMICS, Ecole Nationale des Ponts et
  Chauss\'ees,\\6/8 Av. Blaise Pascal, 77455 Marne-la-Vall\'ee Cedex 02,
  France \\
E-mail: pommaret@cermics.enpc.fr , jean-francois.pommaret@wanadoo.fr\\
URL: http://cermics.enpc.fr/$\sim$pommaret/home.html}
\date{  }
\begin{document}
\maketitle

\noindent
{\bf 1  INTRODUCTION}:\\

    Let us revisit briefly the foundation of n-dimensional elasticity theory as it can be found today in any textbook, restricting our study to $n=2$ for simplicity. If $x=(x^1,x^2)$ is a point in the plane and $\xi=({\xi}^1(x),{\xi}^2(x))$ is the displacement vector, lowering the indices by means of the Euclidean metric, we may introduce the "small" deformation tensor $\epsilon = ({\epsilon}_{ij}={\epsilon}_{ij}=(1/2)({\partial}_i{\xi}_j+{\partial}_j{\xi}_i))$ with $n(n+1)/2=3$ (independent) {\it components} $({\epsilon}_{11}, {\epsilon}_{12}={\epsilon}_{21},{\epsilon}_{22})$. If we study a part of a deformed body, for example a thin elastic plane sheet, by means of a variational principle, we may introduce the local density of free energy $\varphi (\epsilon)=\varphi({\epsilon}_{ij}{\mid} i\leq j)=\varphi ({\epsilon}_{11},{\epsilon}_{12},{\epsilon}_{22})$ and vary the total free energy 
$F={\int}\varphi (\epsilon)dx$ with $dx=dx^1\wedge dx^2$ by introducing ${\sigma}^{ij}=\partial \varphi/\partial {\epsilon}_{ij}$ for $i\leq j$ in order to obtain $\delta F={\int}({\sigma}^{11}\delta {\epsilon}_{11}+{\sigma}^{12}\delta {\epsilon}_{12}+{\sigma}^{22}\delta{\epsilon}_{22})dx$. Accordingly, the "decision" to define the stress tensor $\sigma$ by a symmetric matrix with ${\sigma}^{12}={\sigma}^{21}$ is purely artificial within such a variational principle. Indeed, the usual Cauchy device (1828) assumes that each element of a boundary surface is acted on by a surface density of force ${\vec{\sigma}}$ with a linear dependence $\vec{\sigma}=({\sigma}^{ir}(x)n_r)$ on the outward normal unit vector $\vec{n}=(n_r)$ and does not make any assumption on the stress tensor. It is only by an equilibrium of forces and couples, namely the well known {\it phenomenological static torsor equilibrium}, that one can "prove" the symmetry of $\sigma$. However, even if we assume this symmetry, {\it we now need the different summation} ${\sigma}^{ij}\delta{\epsilon}_{ij}= {\sigma}^{11}\delta{\epsilon}_{11}+2{\sigma}^{12}\delta{\epsilon}_{12}+{\sigma}^{22}\delta{\epsilon}_{22}={\sigma}^{ir}{\partial}_r\delta {\xi}_i$. An integration by part and a change of sign produce the volume integral ${\int}{\partial}_r({\sigma}^{ir})\delta {\xi}_idx$ leading to the stress equations ${\partial}_r{\sigma}^{ir}=0$. {\it The classical approach to elasticity theory, based on invariant theory with respect to the group of rigid motions, cannot therefore describe equilibrium of torsors by means of a variational principle where the proper torsor concept is totally lacking}.\\

     There is another equivalent procedure dealing with a {\it variational calculus with constraint}. Indeed, as we shall see in Section 7, the deformation tensor is not any symmetric tensor as it must satisfy $n^2(n^2-1)/12$ compatibility conditions (CC), that is only ${\partial}_{22}{\epsilon}_{11}+{\partial}_{11}{\epsilon}_{22}-2{\partial}_{12}{\epsilon}_{12}=0$ when $n=2$. In this case, introducing the {\it Lagrange multiplier } $-\phi$ for convenience, {\it we have to vary} $\int({\varphi}(\epsilon)-\phi({\partial}_{22}{\epsilon}_{11}+{\partial}_{11}{\epsilon}_{22}-2{\partial}_{12}{\epsilon}_{12}))dx$ {\it for an arbitrary} $\epsilon$. A double integration by parts now provides the parametrization ${\sigma}^{11}={\partial}_{22}\phi,{\sigma}^{12}={\sigma}^{21}=-{\partial}_{12}\phi, {\sigma}^{22}={\partial}_{11}\phi$ of the stress equations by means of the Airy function $\phi$ and the {\it formal adjoint} of the CC, {\it on the condition to observe that we have in fact} $2{\sigma}^{12}=-2{\partial}_{12}\phi$ as another way to understand the deep meaning of the factor "2" in the summation. In arbitrary dimension, it just remains to notice that the above compatibility conditions are nothing else but the linearized Riemann tensor in Riemannan geometry, a crucial mathematical tool in the theory of general relativity.\\

    It follows that the only possibility to revisit the foundations of engineering and mathematical physics is to use new mathematical methods, namely the theory of systems of partial differential equations and Lie pseudogroups developped by D.C. Spencer and coworkers during the period 1960-1975. In particular, Spencer invented the first order operator now wearing his name in order to bring in a canonical way the formal study of systems of ordinary differential (OD) or partial differential (PD) equations to that of equivalent first order systems. However, despite its importance, the {\it Spencer operator} is rarely used in mathematics today and, up to our knowledge, has never been used in engineering or mathematical physics. The main reason for such a situation is that the existing papers, largely based on hand-written lecture notes given by Spencer to his students (the author was among them in 1969) are quite technical and the problem also lies in the only "accessible" book "Lie equations" he published in 1972 with A. Kumpera. Indeed, the reader can easily check by himself that {\it the core of this book has nothing to do with its introduction} recalling known differential geometric concepts on which most of physics is based today. \\

     The first and technical purpose of this chapter, an extended version of a lecture at the second workshop on Differential Equations by Algebraic Methods (DEAM2, february 9-11, 2011, Linz, Austria), is to recall briefly its definition, both in the framework of systems of linear ordinary or partial differential equations and in the framework of differential modules. The local theory of Lie pseudogroups and the corresponding non-linear framework are also presented for the first time in a rather elementary manner though it is a difficult task.\\
     
     The second and central purpose is to prove that the use of the Spencer operator constitutes the {\it common secret} of the three following famous books published about at the same time in the beginning of the last century, though they do not seem to have anything in common at first sight as they are successively dealing with the foundations of elasticity theory, commutative algebra, electromagnetism (EM) and general relativity (GR):\\  
    
\noindent
[C] E. and F. COSSERAT: "Th\'eorie des Corps D\'eformables", Hermann, Paris, 1909.\\
\noindent
[M] F.S. MACAULAY: "The Algebraic Theory of Modular Systems", Cambridge, 1916.\\
\noindent
[W] H. WEYL: "Space, Time, Matter", Springer, Berlin, 1918 (1922, 1958; Dover, 1952).\\

   Meanwhile we shall point out the striking importance of the second book for studying {\it identifiability} in control theory. We shall also obtain from the previous results the group theoretical unification of finite elements in engineering sciences (elasticity, heat, electromagnetism), solving the {\it torsor problem} and recovering in a purely mathematical way known {\it field-matter coupling phenomena} (piezzoelectricity, photoelasticity, streaming birefringence, viscosity, ...).\\
   
   As a byproduct and though disturbing it may be, the third and perhaps essential purpose is to prove that {\it these unavoidable new differential and homological methods contradict the existing mathematical foundations of both engineering (continuum mechanics, electromagnetism) and mathematical (gauge theory, general relativity) physics}.\\
   
    Many explicit examples will illustate this chapter which is deliberately written in a rather self-contained way to be accessible to a large audience, which does not mean that it is elementary in view of the number of new concepts that must be patched together. However, the reader must never forget that {\it each formula} appearing in this new general framework has been used explicitly or implicitly in [C], [M] and [W] for a mechanical, mathematical or physical purpose. \\

\noindent
{\bf 2   FROM LIE GROUPS TO LIE PSEUDOGROUPS}\\

   Evariste Galois (1811-1832) introduced the word "{\it group}" for the first time in 1830. Then the group concept slowly passed from algebra (groups of permutations) to geometry (groups of transformations). It is only in 1880 that Sophus Lie (1842-1899) studied the groups of transformations depending on a finite number of parameters and now called {\it Lie groups of transformations}. Let $X$ be a manifold with local coordinates $x=(x^1, ... , x^n)$ and $G$ be a {\it Lie group}, that is another manifold with local coordinates $a=(a^1, ... , a^p)$ called {\it parameters} with a {\it composition} $G\times G \rightarrow G: (a,b)\rightarrow ab$, an {\it inverse} $G \rightarrow G: a \rightarrow a^{-1}$ and an {\it identity} $e\in G$ satisfying:\\
\[(ab)c=a(bc)=abc,\hspace{1cm} aa^{-1}=a^{-1}a=e,\hspace{1cm} ae=ea=a,\hspace{1cm} \forall a,b,c \in G \]

\noindent
{\bf DEFINITION} 2.1: $G$ is said to {\it act} on $X$ if there is a map $X\times G \rightarrow X: (x,a) \rightarrow y=ax=f(x,a)$ such that $(ab)x=a(bx)=abx, \forall a,b\in G, \forall x\in X$ and, for simplifying the notations, we shall use global notations even if only local actions are existing. The set $G_x=\{a\in G\mid ax=x\}$ is called the {\it isotropy subgroup} of $G$ at $x\in X$. The action is said to be {\it effective} if $ax=x, \forall x\in X\Rightarrow a=e$. A subset $S\subset X$ is said to be {\it invariant} under the action of $G$ if $aS\subset S,\forall a\in G$ and the {\it orbit} of $x\in X$ is the invariant subset $Gx=\{ax\mid a\in G\}\subset X$. If $G$ acts on two manifolds $X$ and $Y$, a map $f:X\rightarrow Y$ is said to be {\it equivariant} if $f(ax)=af(x), \forall x\in X, \forall a\in G$.\\

  For reasons that will become clear later on, it is often convenient to introduce the {\it graph} $X\times G\rightarrow X\times X: (x,a)\rightarrow (x,y=ax)$ of the action. In the product $X\times X$, the first factor is called the {\it source} while the second factor is called the {\it target}. \\
  
\noindent
{\bf DEFINITION} 2.2: The action is said to be {\it free} if the graph is injective and {\it transitive} if the graph is surjective. The action is said to be {\it simply transitive} if the graph is an isomorphism and $X$ is said to be a {\it principal homogeneous space} (PHS) for $G$. \\

In order to fix the notations, we quote without any proof the "{\it Three Fundamental Theorems of Lie}" that will be of constant use in the sequel ([26]):\\

\noindent
{\bf FIRST FUNDAMENTAL THEOREM} 2.3: The orbits $x=f(x_0,a)$ satisfy the system of PD equations $\partial x^i/\partial a^{\sigma}= {\theta}^i_{\rho}(x){\omega}^{\rho}_{\sigma}(a)$  with $det(\omega)\neq 0$. The vector fields ${\theta}_{\rho}={\theta}^i_{\rho}(x){\partial}_i$ are called {\it infinitesimal generators} of the action and are linearly independent over the constants when the action is effective.\\

If $X$ is a manifold, we denote as usual by $T=T(X)$ the {\it tangent bundle} of $X$, by $T^*=T^*(X)$ the {\it cotangent bundle}, by ${\wedge}^rT^*$ the {\it bundle of r-forms} and by $S_qT^*$ the {\it bundle of q-symmetric tensors}. More generally, let $\cal{E}$ be a {\it fibered manifold}, that is a manifold with local coordinates $(x^i,y^k)$ for $i=1,...,n$ and $k=1,...,m$ simply denoted by $(x,y)$, {\it projection} $\pi:{\cal{E}}\rightarrow X:(x,y)\rightarrow (x)$ and changes of local coordinates $\bar{x}=\varphi(x), \bar{y}=\psi(x,y)$. If $\cal{E}$ and $\cal{F}$ are two fibered manifolds over $X$ with respective local coordinates $(x,y)$ and $(x,z)$, we denote by ${\cal{E}}{\times}_X{\cal{F}}$ the {\it fibered product} of $\cal{E}$ and $\cal{F}$ over $X$ as the new fibered manifold over $X$ with local coordinates $(x,y,z)$. We denote by $f:X\rightarrow {\cal{E}}: (x)\rightarrow (x,y=f(x))$ a global {\it section} of $\cal{E}$, that is a map such that $\pi\circ f=id_X$ but local sections over an open set $U\subset X$ may also be considered when needed. Under a change of coordinates, a section transforms like $\bar{f}(\varphi(x))=\psi(x,f(x))$ and the derivatives transform like:\\
\[   \frac{\partial{\bar{f}}^l}{\partial{\bar{x}}^r}(\varphi(x)){\partial}_i{\varphi}^r(x)=\frac{\partial{\psi}^l}{\partial x^i}(x,f(x))+\frac{\partial {\psi}^l}{\partial y^k}(x,f(x)){\partial}_if^k(x)  \]
We may introduce new coordinates $(x^i,y^k,y^k_i)$ transforming like:\\
\[ {\bar{y}}^l_r{\partial}_i{\varphi}^r(x)=\frac{\partial{\psi}^l}{\partial x^i}(x,y)+\frac{\partial {\psi}^l}{\partial y^k}(x,y)y^k_i  \]
We shall denote by $J_q({\cal{E}})$ the {\it q-jet bundle} of $\cal{E}$ with local coordinates $(x^i, y^k, y^k_i, y^k_{ij},...)=(x,y_q)$ called {\it jet coordinates} and sections $f_q:(x)\rightarrow (x,f^k(x), f^k_i(x), f^k_{ij}(x), ...)=(x,f_q(x))$ transforming like the sections $j_q(f):(x) \rightarrow (x,f^k(x), {\partial}_if^k(x), {\partial}_{ij}f^k(x), ...)=(x,j_q(f)(x))$ where both $f_q$ and $j_q(f)$ are over the section $f$ of $\cal{E}$. Of course $J_q({\cal{E}})$ is a fibered manifold over $X$ with projection ${\pi}_q$ while $J_{q+r}({\cal{E}})$ is a fibered manifold over $J_q({\cal{E}})$ with projection ${\pi}^{q+r}_q, \forall r\geq 0$.\\

\noindent
{\bf DEFINITION} 2.4: A {\it system} of order $q$ on $\cal{E}$ is a fibered submanifold ${\cal{R}}_q\subset J_q({\cal{E}})$ and a {\it solution} of ${\cal{R}}_q$ is a section $f$ of $\cal{E}$ such that $j_q(f)$ is a section of ${\cal{R}}_q$.\\

\noindent
{\bf DEFINITION} 2. 5: When the changes of coordinates have the linear form $\bar{x}=\varphi(x),\bar{y}= A(x)y$, we say that $\cal{E}$ is a {\it vector bundle} over $X$ and denote for simplicity a vector bundle and its set of sections by the same capital letter $E$. When the changes of coordinates have the form $\bar{x}=\varphi(x),\bar{y}=A(x)y+B(x)$ we say that $\cal{E}$ is an {\it affine bundle} over $X$ and we define the {\it associated vector bundle} $E$ over $X$ by the local coordinates $(x,v)$ changing like $\bar{x}=\varphi(x),\bar{v}=A(x)v$.\\

\noindent
{\bf DEFINITION} 2.6: If the tangent bundle $T({\cal{E}})$ has local coordinates $(x,y,u,v)$ changing like ${\bar{u}}^j={\partial}_i{\varphi}^j(x)u^i, {\bar{v}}^l=\frac{\partial {\psi}^l}{\partial x^i}(x,y)u^i+\frac{\partial {\psi}^l}{\partial y^k}(x,y)v^k$, we may introduce the {\it vertical bundle} $V({\cal{E}})\subset T({\cal{E}})$ as a vector bundle over $\cal{E}$ with local coordinates $(x,y,v)$ obtained by setting $u=0$ and changes ${\bar{v}}^l=\frac{\partial {\psi}^l}{\partial y^k}(x,y)v^k$. Of course, when $\cal{E}$ is an affine bundle with associated vector bundle $E$ over $X$, we have $V({\cal{E}})={\cal{E}}\times_XE$.\\

For a later use, if $\cal{E}$ is a fibered manifold over $X$ and $f$ is a section of $\cal{E}$, we denote by $f^{-1}(V({\cal{E}}))$ the {\it reciprocal image} of $V({\cal{E}})$ by $f$ as the vector bundle over $X$ obtained when replacing $(x,y,v)$ by $(x,f(x),v) $ in each chart. It is important to notice in variational calculus that a {\it variation} $\delta f$ of $f$ is such that $\delta f(x)$, as a vertical vector field not necessary "small", is a section of this vector bundle and that $(f,\delta f)$ is nothing else than a section of $V({\cal{E}})$ over $X$.

We now recall a few basic geometric concepts that will be constantly used. First of all, if $\xi,\eta\in T$, we define their {\it bracket} $[\xi,\eta]\in T$ by the local formula $([\xi,\eta])^i(x)={\xi}^r(x){\partial}_r{\eta}^i(x)-{\eta}^s(x){\partial}_s{\xi}^i(x)$ leading to the {\it Jacobi identity} $[\xi,[\eta,\zeta]]+[\eta,[\zeta,\xi]]+[\zeta,[\xi,\eta]]=0, \forall \xi,\eta,\zeta \in T$ allowing to define a {\it Lie algebra} and to the useful formula $[T(f)(\xi),T(f)(\eta)]=T(f)([\xi,\eta])$ where $T(f):T(X)\rightarrow T(Y)$ is the tangent mapping of a map $f:X\rightarrow Y$.\\

\noindent
{\bf SECOND FUNDAMENTAL THEOREM} 2.7: If ${\theta}_1,...,{\theta}_p$ are the infinitesimal generators of the effective action of a lie group $G$ on $X$, then $[{\theta}_{\rho},{\theta}_{\sigma}]=c^{\tau}_{\rho\sigma}{\theta}_{\tau}$ where the $c^{\tau}_{\rho\sigma}$ are the {\it structure constants} of a Lie algebra of vector fields which can be identified with ${\cal{G}}=T_e(G)$.\\

When $I=\{ i_1< ... < i_r\}$ is a multi-index, we may set $dx^I=dx^{i_1}\wedge ... \wedge dx^{i_r}$ for describing ${\wedge}^rT^*$ and introduce the {\it exterior derivative} $d:{\wedge}^rT^*\rightarrow {\wedge}^{r+1}T^*:\omega={\omega}_Idx^I \rightarrow d\omega={\partial}_i{\omega}_Idx^i\wedge dx^I$ with $d^2=d\circ d\equiv 0$ in the {\it Poincar\'{e} sequence}:\\
\[  {\wedge}^0T^* \stackrel{d}{\longrightarrow} {\wedge}^1T^* \stackrel{d}{\longrightarrow} {\wedge}^2T^* \stackrel{d}{\longrightarrow} ... \stackrel{d}{\longrightarrow} {\wedge}^nT^* \longrightarrow 0  \]

The {\it Lie derivative} of an $r$-form with respect to a vector field $\xi\in T$ is the linear first order operator ${\cal{L}}(\xi)$ linearly depending on $j_1(\xi)$ and uniquely defined by the following three properties:\\
1) ${\cal{L}}(\xi)f=\xi.f={\xi}^i{\partial}_if, \forall f\in {\wedge}^0T^*=C^{\infty}(X)$.\\
2) ${\cal{L}}(\xi)d=d{\cal{L}}(\xi)$.\\
3) ${\cal{L}}(\xi)(\alpha\wedge \beta)=({\cal{L}}(\xi)\alpha)\wedge \beta+\alpha\wedge ({\cal{L}}(\xi) \beta), \forall \alpha,\beta \in \wedge T^*$.\\
It can be proved that ${\cal{L}}(\xi)=i(\xi)d+di(\xi)$ where $i(\xi)$ is the {\it interior multiplication} $(i(\xi)\omega)_{i_1...i_r}={\xi}^i{\omega}_{ii_1...i_r}$ and that $[{\cal{L}}(\xi),{\cal{L}}(\eta)]={\cal{L}}(\xi)\circ {\cal{L}}(\eta)-{\cal{L}}(\eta)\circ {\cal{L}}(\xi)={\cal{L}}([\xi,\eta]), \forall \xi,\eta\in T$.\\

Using crossed-derivatives in the PD equations of the First Fundamental Theorem and introducing the family of $1$-forms ${\omega}^{\tau}={\omega}^{\tau}_{\sigma}(a)da^{\sigma}$ both with the matrix $\alpha={\omega}^{-1}$ of right invariant vector fields, we obtain the {\it compatibility conditions} (CC) expressed by the following corollary where one must care about the sign used:\\

{\bf COROLLARY} 2.8: One has the {\it Maurer-Cartan} (MC) {\it equations} $d{\omega}^{\tau}+c^{\tau}_{\rho\sigma}{\omega}^{\rho}\wedge{\omega}^{\sigma}=0$ or the equivalent relations $[{\alpha}_{\rho},{\alpha}_{\sigma}]=c^{\tau}_{\rho\sigma}{\alpha}_{\tau}$.\\

Applying $d$ to the MC equations and substituting, we obtain the {\it integrability conditions} (IC):\\

\noindent
{\bf THIRD FUNDAMENTAL THEOREM} 2.9: For any Lie algebra $\cal{G}$ defined by structure constants satisfying :\\
\[  c^{\tau}_{\rho\sigma}+c^{\tau}_{\sigma \rho}=0, \hspace{1cm} c^{\lambda}_{\mu\rho}c^{\mu}_{\sigma\tau}+c^{\lambda}_{\mu\sigma}c^{\mu}_{\tau\rho}+c^{\lambda}_{\mu\tau}c^{\mu}_{\rho\sigma}=0  \]
one can construct an analytic group $G$ such that ${\cal{G}}=T_e(G)$.\\

\noindent
{\bf EXAMPLE} 2.10: Considering the affine group of transformations of the real line $y=a^1x+a^2$, we obtain ${\theta}_1=x{\partial}_x, {\theta}_2={\partial}_x \Rightarrow [{\theta}_1,{\theta}_2]=-{\theta}_2$ and thus ${\omega}^1=(1/{a^1})da^1, {\omega}^2=da^2-(a^2/{a^1})da^1\Rightarrow d{\omega}^1=0, d{\omega}^2-{\omega}^1\wedge{\omega}^2=0$.\\

Only ten years later Lie discovered that the Lie groups of transformations are only particular examples of a wider class of groups of transformations along the following definition where $aut(X)$ denotes the group of all local diffeomorphisms of $X$:\\

\noindent
{\bf DEFINITION} 2.11: A {\it Lie pseudogroup} of transformations $\Gamma\subset aut(X)$ is a group of transformations solutions of a system of OD or PD equations such that, if $y=f(x)$ and $z=g(y)$ are two solutions, called {\it finite transformations}, that can be composed, then $z=g\circ f(x)=h(x)$ and 
$x=f^{-1}(y)=g(y)$ are also solutions while $y=x$ is a solution.\\

  The underlying system may be nonlinear and of high order as we shall see later on. We shall speak of an {\it algebraic pseudogroup} when the system is defined by {\it differential polynomials} that is polynomials in the derivatives. In the case of Lie groups of transformations the system is obtained by differentiating the action law $y=f(x,a)$ with respect to $x$ as many times as necessary in order to eliminate the parameters. Looking for transformations "close" to the identity, that is setting $y=x+t\xi(x)+...$ when $t\ll 1$ is a small constant parameter and passing to the limit $t\rightarrow 0$, we may linearize the above nonlinear {\it system of finite Lie equations} in order to obtain a linear {\it system of infinitesimal Lie equations} of the same order for vector fields. Such a system has the property that, if $\xi,\eta$ are two solutions, then $[\xi,\eta]$ is also a solution. Accordingly, the set $\Theta\subset T$ of solutions of this new system satifies $[\Theta,\Theta]\subset \Theta$ and can therefore be considered as the Lie algebra of $\Gamma$.\\

  Though the collected works of Lie have been published by his student F. Engel at the end of the $19^{th}$ century, these ideas did not attract a large audience because the fashion in Europe was analysis. Accordingly, at the beginning of the $20^{th}$ century and for more than fifty years, only two frenchmen tried to go further in the direction pioneered by Lie, namely Elie Cartan (1869-1951) who is quite famous today and Ernest Vessiot (1865-1952) who is almost ignored today, each one deliberately ignoring the other during his life for a precise reason that we now explain with details as it will surprisingly constitute the heart of this chapter. (The author is indebted to Prof. Maurice Janet (1888-1983), who was a personal friend of Vessiot, for the many documents he gave him, partly published in [25]). Roughly, the idea of many people at that time was to extend the work of Galois along the following scheme of increasing difficulty:\\
1) {\it Galois theory} ([34]): Algebraic equations and permutation groups.\\
2) {\it Picard-Vessiot theory} ([17]): OD equations and Lie groups.\\
3) {\it Differential Galois theory} ([24],[37]): PD equations and Lie pseudogroups.\\
In 1898 Jules Drach (1871-1941) got and published a thesis ([9]) with a jury made by Gaston Darboux, Emile Picard and Henri Poincar\'{e}, the best leading mathematicians of that time. However, despite the fact that it contains ideas quite in advance on his time such as the concept of a "differential field" only introduced by J.-F. Ritt in 1930 ([31]), the jury did not notice that the main central result was wrong: Cartan found the counterexamples, Vessiot recognized the mistake, Paul Painlev\'{e} told it to Picard who agreed but Drach never wanted to acknowledge this fact and was supported by the influent Emile Borel. As a byproduct, everybody flew out of this "affair", never touching again the Galois theory. After publishing a prize-winning paper in 1904 where he discovered for the first time that the differential Galois theory must be a theory of (irreducible) PHS for algebraic pseudogroups, Vessiot remained alone, trying during thirty years to correct the mistake of Drach that we finally corrected in 1983 ([24]).\\

\noindent
{\bf 3   CARTAN VERSUS VESSIOT : THE STRUCTURE EQUATIONS}\\

  We study first the work of Cartan which can be divided into two parts. The first part, for which he invented exterior calculus, may be considered as a tentative to extend the MC equations from Lie groups to Lie pseudogroups. The idea for that is to consider the system of order $q$ and its {\it prolongations} obtained by differentiating the equations as a way to know certain derivatives called {\it principal} from all the other arbitrary ones called {\it parametric} in the sense of Janet ([13]). Replacing the derivatives by jet coordinates, we may try to copy the procedure leading to the MC equations by using a kind of "composition" and "inverse" on the jet coordinates. For example, coming back to the last definition, we get successively:\\
  \[  \frac{\partial h}{\partial x}=\frac{\partial g}{\partial y}\frac{\partial f}{\partial x}, \hspace{5mm}\frac{{\partial}^2h}{\partial x^2}=\frac{{\partial}^2g}{\partial y^2}\frac{\partial f}{\partial x}\frac{\partial f}{\partial x}+\frac{\partial g}{\partial y}\frac{{\partial}^2f}{\partial x^2}, ...   \]
Now if $g=f^{-1}$ then $g\circ f=id$ and thus $\frac{\partial g}{\partial y}\frac{\partial f}{\partial x}=1,...$ while the new identity $id_q=j_q(id)$ is made by the successive derivatives of $y=x$, namely $(1,0,0,...)$. These {\it awfully complicated computations} bring a lot of structure constants and have been so much superseded by the work of Donald C. Spencer (1912-2001) ([11],[12],[18],[33]) that, in our opinion based on thirty years of explicit computations, this tentative has only been used for classification problems and is not useful for applications compared to the results of the next sections. In a single concluding sentence, {\it Cartan has not been able to "go down" to the base manifold} $X$ {\it while Spencer did succeed fifty years later}.\\

We shall now describe the second part with more details as it has been (and still is !) the crucial tool used in both engineering (analytical and continuum mechanics) and mathematical (gauge theory and general relativity) physics in an absolutely contradictory manner. We shall try to use the least amount of mathematics in order to prepare the reader for the results presented in the next sections. For this let us start with an elementary experiment that will link at once continuum mechanics and gauge theory in an unusual way. Let us put a thin elastic rectilinear rubber band along the $x$ axis on a flat table and translate it along itself. The band will remain identical as no deformation can be produced by this constant translation. However, if we move each point continuously along the same direction but in a point depending way, for example fixing one end and pulling on the other end, there will be of course a deformation of the elastic band according to the Hooke law. It remains to notice that a constant translation can be written in the form $y=x+a^2$ as in Example 2.10 while a point depending translation can be written in the form $y=x+a^2(x)$ which is written in any textbook of continuum mechanics in the form $y=x+\xi(x)$ by introducing the {\it displacement vector} $\xi$. However nobody could even imagine to make $a^1$ also point depending and to consider $y=a^1(x)x+a^2(x)$ as we shall do in Example 7.9.We also notice that the movement of a rigid body in space may be written in the form $y=a(t)x+b(t)$ where now $a(t)$ is a time depending orthogonal matrix and $b(t)$ is a time depending vector. What makes all the difference between the two examples is that the group is {\it acting} on $x$ in the first but {\it not acting} on $t$ in the second. Finally, a point depending rotation or dilatation is not accessible to intuition and the general theory must be done in the following manner.\\

If $X$ is a manifold and $G$ is a lie group {\it not acting necessarily} on $X$, let us consider maps $a:X\rightarrow G: (x)\rightarrow (a(x))$ or equivalently sections of the trivial (principal) bundle $X\times G$ over $X$. If $x+dx$ is a point of $X$ close to $x$, then $T(a)$ will provide a point $a+da=a+\frac{\partial a}{\partial x}dx$ close to $a$ on $G$. We may bring $a$ back to $e$ on $G$ by acting on $a$ with $a^{-1}$, {\it  either on the left or on the right}, getting therefore a $1$-form $a^{-1}da=A$ or $daa^{-1}=B$. As $aa^{-1}=e$ we also get $daa^{-1}=-ada^{-1}=-b^{-1}db$ if we set $b=a^{-1}$ as a way to link $A$ with $B$. When there is an action $y=ax$, we have $x=a^{-1}y=by$ and thus $dy=dax=daa^{-1}y$, a result leading through the First Fundamental Theorem of Lie to the equivalent formulas:\\
\[   a^{-1}da=A=({A}^{\tau}_i(x)dx^i=-{\omega}^{\tau}_{\sigma}(b(x)){\partial}_ib^{\sigma}(x)dx^i)  \]
\[   daa^{-1}=B=({B}^{\tau}_i(x)dx^i={\omega}^{\tau}_{\sigma}(a(x)){\partial}_ia^{\sigma}(x)dx^i)  \]
Introducing the induced bracket $[A,A](\xi,\eta)=[A(\xi),A(\eta)]\in {\cal{G}}, \forall \xi,\eta\in T$, we may define the $2$-form $dA-[A,A]=F\in {\wedge}^2T^*\otimes {\cal{G}}$ by the local formula (care to the sign):\\
\[     {\partial}_iA^{\tau}_j(x)-{\partial}_jA^{\tau}_i(x)-c^{\tau}_{\rho\sigma}A^{\rho}_i(x)A^{\sigma}_j(x)=F^{\tau}_{ij}(x)  \]
and obtain from the second fundamental theorem:\\

\noindent
{\bf THEOREM} 3.1: There is a {\it nonlinear gauge sequence}:\\
\[  \begin{array}{ccccc}
X\times G & \longrightarrow & T^*\otimes {\cal{G}} &\stackrel{MC}{ \longrightarrow} & {\wedge}^2T^*\otimes {\cal{G}}  \\
a                & \longrightarrow  &    a^{-1}da=A         &    \longrightarrow & dA-[A,A]=F
\end{array}   \]

Choosing $a$ "close" to $e$, that is $a(x)=e+t\lambda(x)+...$ and linearizing as usual, we obtain the linear operator $d:{\wedge}^0T^*\otimes {\cal{G}}\rightarrow {\wedge}^1T^*\otimes {\cal{G}}:({\lambda}^{\tau}(x))\rightarrow ({\partial}_i{\lambda}^{\tau}(x))$ leading to:\\

\noindent
{\bf COROLLARY} 3.2: There is a {\it linear gauge sequence}:\\ 
\[  {\wedge}^0T^*\otimes {\cal{G}}\stackrel{d}{\longrightarrow} {\wedge}^1T^*\otimes {\cal{G}} \stackrel{d}{\longrightarrow} {\wedge}^2T^*\otimes{\cal{G}} \stackrel{d}{\longrightarrow} ... \stackrel{d}{\longrightarrow} {\wedge}^nT^*\otimes {\cal{G}}\longrightarrow  0   \]
which is the tensor product by $\cal{G}$ of the Poincar\'{e} sequence:\\

\noindent
{\bf REMARK} 3.3: When the physicists C.N. Yang and R.L. Mills created (non-abelian) gauge theory in 1954 ([38],[39]), their work was based on these results which were the only ones known at that time, the best mathematical reference being the well known book by S. Kobayashi and K. Nomizu ([15]). It follows that {\it the only possibility} to describe elecromagnetism (EM) within this framework was to call $A$ the {\it Yang-Mills potential} and $F$ the {\it Yang-Mills field} (a reason for choosing such notations) on the condition to have $dim(G)=1$ in the abelian situation $c=0$ and to use a Lagrangian depending on $F=dA-[A,A]$ in the general case. Accordingly the idea was to select the {\it unitary group} $U(1)$, namely the unit circle in the complex plane with Lie algebra the tangent line to this circle at the unity $(1,0)$. It is however important to notice that the resulting Maxwell equations $dF=0$ have no equivalent in the non-abelian case $c\neq 0$.\\

  Just before Albert Einstein visited Paris in 1922, Cartan published many short Notes ([5]) announcing long papers ([6]) where he selected $G$ to be the Lie group involved in the Poincar\'{e} (conformal) group of space-time preserving (up to a function factor) the Minkowski metric $\omega=(dx^1)^2+(dx^2)^2+(dx^3)^2-(dx^4)^2$ with $x^4=ct$ where $c$ is the speed of light. In the first case $F$ is decomposed into two parts, the {\it torsion} as a $2$-form with value in translations on one side and the {\it curvature} as a $2$-form with value in rotations on the other side. This result was looking coherent {\it at first sight} with the Hilbert variational scheme of general relativity (GR) introduced by Einstein in 1915 ([21],[38]) and leading to a Lagrangian depending on $F=dA-[A,A]$ as in the last remark.\\
  
  In the meantime, Poincar\'{e} developped an invariant variational calculus ([22]) which has been used again without any quotation, successively by G. Birkhoff and V. Arnold (compare [4], 205-216 with [2], 326, Th 2.1).  A particular case is well known by any student in the analytical mechanics of rigid bodies. Indeed, using standard notations, the movement of a rigid body is described in a fixed Cartesian frame by the formula $x(t)=a(t)x_0+b(t)$ where $a(t)$ is a $3\times 3$ time dependent orthogonal matrix (rotation) and $b(t)$ a time depending vector (translation) as we already said. Differentiating with respect to time by using a dot, the {\it absolute speed} is $v={\dot{x}}(t)={\dot{a}}(t)x_0+{\dot{b}}(t)$ and we obtain the {\it relative speed} $a^{-1}(t)v=a^{-1}(t){\dot{a}}(t)x_0+a^{-1}(t){\dot{b}}(t)$ by projection in a frame fixed in the body. Having in mind Example 2.10, it must be noticed that the so-called {\it Eulerian speed} $v=v(x,t)=\dot{a}a^{-1}x+\dot{b}-\dot{a}a^{-1}b$ only depends on the $1$-form $B=(\dot{a}a^{-1},\dot{b}-\dot{a}a^{-1}b)$. The Lagrangian (kinetic energy in this case) is thus a quadratic function of the $1$-form $A=(a^{-1}\dot{a},a^{-1}\dot{b})$ where $a^{-1}\dot{a}$ is a $3\times 3$ skew symmetric time depending matrix. Hence, "{\it surprisingly}", this result is not coherent at all with EM where the Lagrangian is the quadratic expression $(\epsilon/2)E^2-(1/2{\mu})B^2$ because the electric field $\vec{E}$ and the magnetic field $\vec{B}$ are combined in the EM field $F$ as a $2$-form satisfying the first set of Maxwell equations $dF=0$. The dielectric constant $\epsilon$ and the magnetic constant $\mu$ are leading to the electric induction $\vec{D}=\epsilon\vec{E}$ and the magnetic induction $\vec{H}=(1/\mu)\vec{B}$ in the second set of Maxwell equations. In view of the existence of well known field-matter couplings such as piezoelectricity and photoelasticity that will be described later on, such a situation is contradictory as it should lead to put on equal footing $1$-forms and $2$-forms contrary to any unifying mathematical scheme but no other substitute could have been provided at that time.\\

Let us now turn to the other way proposed by Vessiot in 1903 ([36]) and 1904 ([37]). Our purpose is only to sketch the main results that we have obtained in many books ([23-26], we do not know other references) and to illustrate them by a series of specific examples, asking the reader to imagine any link with what has been said.\\

1) If ${\cal{E}}=X\times X$, we shall denote by ${\Pi}_q={\Pi}_q(X,X)$ the open subfibered manifold of $J_q(X\times X)$ defined independently of the coordinate system by $det(y^k_i)\neq 0$ with {\it source projection} ${\alpha}_q:{\Pi}_q\rightarrow X:(x,y_q)\rightarrow (x)$ and {\it target projection} ${\beta}_q:{\Pi}_q\rightarrow X:(x,y_q)\rightarrow (y)$. We shall sometimes introduce a copy $Y$ of $X$ with local coordinates $(y)$ in order to avoid any confusion between the source and the target manifolds. Let us start with a Lie pseudogroup $\Gamma\subset aut(X)$ defined by a system ${\cal{R}}_q\subset {\Pi}_q$ of order $q$. In all the sequel we shall suppose that the system is involutive (see next section) and that $\Gamma$ is {\it transitive} that is $\forall x,y\in X, \exists f\in \Gamma, y=f(x)$ or, equivalently, the map $({\alpha}_q,{\beta}_q):{\cal{R}}_q\rightarrow X\times X:(x,y_q)\rightarrow (x,y)$ is surjective.\\

2) The lie algebra $\Theta\subset T$ of infinitesimal transformations is then obtained by linearization, setting $y=x+t\xi(x)+...$ and passing to the limit $t\rightarrow 0$ in order to obtain the linear involutive system $R_q= id^{-1}_q(V({\cal{R}}_q))\subset J_q(T)$ by reciprocal image with $\Theta=\{\xi \in T{\mid}j_q(\xi) \in R_q\}$.  \\

3) Passing from source to target, we may {\it prolong} the vertical infinitesimal transformations $\eta={\eta}^k(y)\frac{\partial}{\partial y^k}$ to the jet coordinates up to order $q$ in order to obtain:\\
\[   {\eta}^k(y)\frac{\partial}{\partial y^k}+\frac{\partial {\eta}^k}{\partial y^r}y^r_i\frac{\partial}{\partial y^k_i}+(\frac{{\partial}^2{\eta}^k}{\partial y^r\partial y^s}y^r_iy^s_j+\frac{\partial {\eta}^k}{\partial y^r}y^r_{ij})\frac{\partial}{\partial y^k_{ij}}+...    \]
where we have replaced $j_q(f)(x)$ by $y_q$, each component beeing the "formal" derivative of the previous one .\\

4) As $[\Theta,\Theta]\subset \Theta$, we may use the Frobenius theorem in order to find a generating fundamental set of {\it differential invariants} $\{{\Phi}^{\tau}(y_q)\}$ up to order $q$ which are such that ${\Phi}^{\tau}({\bar{y}}_q)={\Phi}^{\tau}(y_q)$ by using the chain rule for derivatives whenever $\bar{y}=g(y)\in \Gamma$ acting now on $Y$. Of course, in actual practice {\it one must use sections of} $R_q$ {\it instead of solutions} but it is only in section 6 that we shall see why the use of the Spencer operator will be crucial for this purpose. Specializing the ${\Phi}^{\tau}$ at $id_q(x)$ we obtain the {\it Lie form} ${\Phi}^{\tau}(y_q)={\omega}^{\tau}(x)$ of ${\cal{R}}_q$.\\

5) The main discovery of Vessiot, fifty years in advance, has been to notice that the prolongation at order $q$  of any horizontal vector field $\xi={\xi}^i(x)\frac{\partial}{\partial x^i}$ commutes with the prolongation at order $q$ of any vertical vector field $\eta={\eta}^k(y)\frac{\partial}{\partial y^k}$, exchanging therefore the differential invariants. Keeping in mind the well known property of the Jacobian determinant while passing to the finite point of view, any (local) transformation $y=f(x)$ can be lifted to a (local) transformation of the differential invariants between themselves of the form $u\rightarrow \lambda(u,j_q(f)(x))$ allowing to introduce a {\it natural bundle} $\cal{F}$ over $X$ by patching changes of coordinates $\bar{x}=\varphi(x), \bar{u}=\lambda(u,j_q(\varphi)(x))$. A section $\omega$ of $\cal{F}$ is called a {\it geometric object} or {\it structure} on $X$ and transforms like ${\bar{\omega}}(f(x))=\lambda(\omega(x),j_q(f)(x))$ or simply $\bar{\omega}=j_q(f)(\omega)$. This is a way to generalize vectors and tensors ($q=1$) or even connections ($q=2$). As a byproduct we have $\Gamma=\{f\in aut(X){\mid} {\Phi}_{\omega}(j_q(f))=j_q(f)^{-1}(\omega)=\omega\}$ as a new way to write out the Lie form and we may say that $\Gamma$ {\it preserves} $\omega$. We also obtain ${\cal{R}}_q=\{f_q\in {\Pi}_q{\mid} f_q^{-1}(\omega)=\omega\}$. Coming back to the infinitesimal point of view and setting $f_t=exp(t\xi)\in aut(X), \forall \xi\in T$, we may define the {\it ordinary Lie derivative} with value in ${\omega}^{-1}(V({\cal{F}}))$ by the formula :\\
\[   {\cal{D}}\xi={\cal{D}}_{\omega}\xi={\cal{L}}(\xi)\omega=\frac{d}{dt}j_q(f_t)^{-1}(\omega){\mid}_{t=0} \Rightarrow \Theta=\{\xi\in T{\mid}{\cal{L}}(\xi)\omega=0\}      \]
while we have $x\rightarrow x+t\xi(x)+...\Rightarrow u^{\tau}\rightarrow u^{\tau}+t{\partial}_{\mu}{\xi}^kL^{\tau\mu}_k(u)+...$ where $\mu=({\mu}_1,...,{\mu}_n)$ is a multi-index as a way to write down the system of infinitesimal Lie equations in the {\it Medolaghi form}:\\
\[     {\Omega}^{\tau}\equiv ({\cal{L}}(\xi)\omega)^{\tau}\equiv -L^{\tau\mu}_k(\omega(x)){\partial}_{\mu}{\xi}^k+{\xi}^r{\partial}_r{\omega}^{\tau}(x)=0    \]

7) By analogy with "special" and "general" relativity, we shall call the given section {\it special} and any other arbitrary section {\it general}. The problem is now to study the formal properties of the linear system just obtained with coefficients only depending on $j_1(\omega)$, exactly like L.P. Eisenhart did for ${\cal{F}}=S_2T^*$ when finding the constant Riemann curvature condition for a metric $\omega$ with 
$det(\omega)\neq 0$ ([26], Example 10, p 249). Indeed, if any expression involving $\omega$ and its derivatives is a scalar object, it must reduce to a constant because $\Gamma$ is assumed to be transitive and thus cannot be defined by any zero order equation. Now one can prove that the CC for $\bar{\omega}$, thus for $\omega$ too, only depend on the $\Phi$ and take the quasi-linear symbolic form $v\equiv I(u_1)\equiv A(u)u_x+B(u)=0$, allowing to define an affine subfibered manifold ${\cal{B}}_1\subset J_1({\cal{F}})$ over $\cal{F}$. 
Now, if one has two sections $\omega$ and $\bar{\omega}$ of $\cal{F}$, the {\it equivalence problem} is to look for $f\in aut(X)$ such that $j_q(f)^{-1}(\omega)=\bar{\omega}$. When the two sections satisfy the same CC, the problem is sometimes locally possible (Lie groups of transformations, Darboux problem in analytical mechanics,...) but sometimes not ([23], p 333).\\

8) Instead of the CC for the equivalence problem, let us look for the {\it integrability conditions} (IC) for the system of infinitesimal Lie equations and suppose that, for the given section, all the equations of order $q+r$ are obtained by differentiating $r$ times only the equations of order $q$, then it was claimed by Vessiot 
([36] with no proof, see [26], p 209) that such a property is held if and only if there is an equivariant section 
$c:{\cal{F}}\rightarrow {\cal{F}}_1:(x,u)\rightarrow (x,u,v=c(u))$ where ${\cal{F}}_1=J_1({\cal{F}})/{\cal{B}}_1$ is a natural vector bundle over $\cal{F}$ with local coordinates $(x,u,v)$. Moreover, any such equivariant section depends on a finite number of constants $c$ called {\it structure constants} and the IC for the {\it Vessiot structure equations} $I(u_1)=c(u)$ are of a polynomial form $J(c)=0$.\\

9) Finally, when $Y$ is no longer a copy of $X$, a system ${\cal{A}}_q\subset J_q(X\times Y)$ is said to be an {\it automorphic system} for a Lie pseudogroup $\Gamma\subset aut(Y)$ if, whenever $ y=f(x)$ and $\bar{y}=\bar{f}(x)$ are two solutions, then there exists one and only one transformation $\bar{y}=g(y)\in \Gamma$ such that $\bar{f}=g\circ f$. Explicit tests for checking such a property formally have been given in [24] and can be implemented on computer in the differential algebraic framework.\\

\noindent
{\bf EXAMPLE} 3.4: ({\it Principal homogeneous structure}) When $\Gamma$ is made by the translations $y^i=x^i+a^i$, the Lie form is ${\Phi}^k_i(y_1)\equiv y^k_i={\delta}^k_i$ (Kronecker symbol) and the linearization is ${\partial}_i{\xi}^k=0$. The natural bundle is ${\cal{F}}=T^*{\times}_X... {\times}_X T^*$ ($n$ times) with $det(\omega)\neq 0$ and the general Medolaghi form is ${\omega}^{\tau}_r{\partial}_i{\xi}^r+{\xi}^r{\partial}_r{\omega}^{\tau}_i=0\Leftrightarrow [\xi,{\alpha}_{\tau}]=0$ with $\tau=1,...,n$ if $\alpha=({\alpha}^i_{\tau})={\omega}^{-1}$. Using crossed drivatives, one finally gets the zero order equations: \\
\[     {\xi}^r{\partial}_r({\alpha}^i_{\rho}(x){\alpha}^j_{\sigma}(x)({\partial}_i{\omega}^{\tau}_j(x)-{\partial}_j{\omega}^{\tau}_i(x)))=0  \]
leading therefore (up to sign) to the $n^2(n-1)/2$ Vessiot structure equations:\\
\[     {\partial}_i{\omega}^{\tau}_j(x)-{\partial}_j{\omega}^{\tau}_i(x)=c^{\tau}_{\rho\sigma}{\omega}^{\rho}_i(x){\omega}^{\sigma}_j(x) \]
This result proves that the MC equations are only examples of the Vessiot structure equations. We finally explain the name given to this structure ([26], p 296). Indeed, when $X$ is a PHS for a lie group $G$, the graph of the action is an isomorphism and we obtain a map $X\times X\rightarrow G: (x,y)\rightarrow (a(x,y))$ leading to a first order system of finite Lie equations $y_x=\frac{\partial f}{\partial x} (x,a(x,y))$. In order to produce a Lie form, let us first notice that the general solution of the system of infinitesimal equations is $\xi={\lambda}^{\tau}{\theta}_{\tau}$ with $\lambda=cst$. Introducing the inverse matrix $(\omega)=({\omega}^{\tau}_i)$ of the {\it reciprocal distribution} $\alpha=\{{\alpha}_{\tau}\}$ made by n vectors commuting with $\{{\theta}_{\tau}\}$, we obtain $\lambda = cst \Leftrightarrow [\xi,\alpha]=0 \Leftrightarrow{\cal{L}}(\xi)\omega=0$.    \\

\noindent
{\bf EXAMPLE} 3.5: ({\it Affine and projective structures of the real line}) In Example 2.10 with $n=1$, the special Lie equations are $\Phi(y_2)\equiv y_{xx}/y_x=0 \Rightarrow {\partial}_{xx}\xi=0$ with $q=2$ and we let the reader check as an exercise that the general Lie equations are:\\
\[  \frac{y_{xx}}{y_x}+\omega(y)y_x=\omega(x)  \Rightarrow  {\partial}_{xx}\xi+\omega(x){\partial}_x\xi+\xi{\partial}_x\omega(x)=0  \]
with no IC. The special section is $\omega(x)=0$. \\
We could study in the same way the group of projective transformations of the real line $y=(ax+b)/(cx+d)$ and get with more work the general lie equations:\\
\[   \frac{y_{xxx}}{y_x}-\frac{3}{2}(\frac{y_{xx}}{y_x})^2+\omega(y)y^2_x=\omega(x) \Rightarrow {\partial}_{xxx}\xi+2\omega(x){\partial}_x\xi+\xi{\partial}_x\omega(x)=0   \]
There is an isomorphism $J_1({\cal{F}}_{aff})\simeq {\cal{F}}_{aff}{\times}_X{\cal{F}}_{proj}: 
j_1(\omega)\rightarrow (\omega,\gamma={\partial}_x\omega-(1/2){\omega}^2)$.\\

\noindent
{\bf EXAMPLE} 3.6: $n=2, q=1, \Gamma=\{y^1=f(x^1), y^2=x^2/(\partial f(x^1)/\partial x^1)\} $ where $f$ is an arbitrary invertible map. The involutive Lie form is:\\
\[   {\Phi}^1(y_1)\equiv y^2y^1_1=x^2, \hspace{1cm} {\Phi}^2(y_1)\equiv y^2y^1_2=0, \hspace{1cm}{\Phi}^3(y_1)\equiv \frac{\partial(y^1,y^2)}{\partial(x^1,x^2)}\equiv y^1_1y^2_2-y^1_2y^2_1=1  \]
We obtain ${\cal{F}}=T^*{\times}_X{\wedge}^2T^*$ and $\omega=(\alpha,\beta)$ where $\alpha$ is a $1$-form and $\beta$ is a $2$-form with special section $\omega=(x^2dx^1, dx^1\wedge dx^2)$. It follows that $d\alpha/\beta$ is a well defined scalar because $\beta\neq 0$. The Vessiot structure equation is $d\alpha=c\beta$ with a single structure constant $c$ which cannot have anything to do with a Lie algebra. Considering the other section $\bar{\omega}=(dx^1,dx^1\wedge dx^2)$, we get $\bar{c}=0$. As $\bar{c}\neq c$, the equivalence problem $j_1(f)^{-1}(\omega)=\bar{\omega}$ cannot even be solved formally.\\

\noindent
{\bf EXAMPLE} 3.7: ({\it Symplectic structure}) With $n=2p, q=1$ and ${\cal{F}}={\wedge}^2T^*$, let $\omega$ be a closed $2$-form of maximum rank, that is $d\omega=0, det(\omega)\neq 0$. The equivalence problem is nothing else than the Darboux problem in analytical mechanics giving the possibility to write locally $\omega=\sum dp\wedge dq$ by using canonical conjugate coordinates $(q,p)=(position, momentum)$.\\

\noindent
{\bf EXAMPLE} 3.8: ({\it Contact structure}) With $n=3, q=1, w=dx^1-x^3dx^2\Rightarrow w\wedge dw=dx^1\wedge dx^2\wedge dx^3$, let us consider $\Gamma=\{f\in aut(X){\mid}j_1(f)^{-1}(w)=\rho w\}$. {\it This is not a Lie form} but we get: \\
\[   j_1(f)^{-1}(dw)=dj_1(f)^{-1}(w)=\rho dw+d\rho\wedge w \Rightarrow j_1(f)^{-1}(w\wedge dw)={\rho}^2 (w\wedge dw)  \]
The corresponding geometric object is thus made by a $1$-form {\it density} $\omega=({\omega}_1,{\omega}_2,{\omega}_3)$ that transforms like a $1$-form up to the division by the square root of the Jacobian determinant. The unusual general Medolaghi form is: \\
\[   {\Omega}_i\equiv {\omega}_r(x){\partial}_i{\xi}^r-(1/2){\omega}_i(x){\partial}_r{\xi}^r+{\xi}^r{\partial}_r{\omega}_i(x)=0\]
In a symbolic way $\omega\wedge d\omega$ is now a scalar and the only Vessiot structure equation is: \\
\[   {\omega}_1({\partial}_2{\omega}_3-{\partial}_3{\omega}_2)+ {\omega}_2({\partial}_3{\omega}_1-{\partial}_1{\omega}_3)+{\omega}_3({\partial}_1{\omega}_2-{\partial}_2{\omega}_1)=c   \]
For the special section $\omega=(1,-x^3,0)$ we have $c=1$. If we choose $\bar{\omega}=(1,0,0)$ we may define $\bar{\Gamma}$ by the system $y^1_2=0, y^1_3=0, y^2_2y^3_3-y^2_3y^3_2=y^1_1$ but now $\bar{c}=0$ and the equivalence problem $j_1(f)^{-1}(\omega)=\bar{\omega}$ cannot even be solved formally. These results can be extended to an arbitrary odd dimension with much more work ([24], p 684).\\

\noindent
{\bf EXAMPLE} 3.9: ({\it Screw and complex structures}) ($n=2, q=1$) In 1878 Clifford introduced abstract numbers of the form $x^1+\epsilon x^2$ with ${\epsilon}^2=0$ in order to study helicoidal movements in the mechanics of rigid bodies. We may try to define functions of these numbers for which a derivative may have a meaning. Thus, if $f(x^1+\epsilon x^2)=f^1(x^1,x^2)+\epsilon f^2(x^1,x^2)$, then we should get:\\
\[ df=(A+\epsilon B)(dx^1+\epsilon dx^2)=Adx^1+\epsilon (Bdx^1+Adx^2)=df^1+\epsilon df^2\]
Accordingly, we have to look for transformations $y^1=f^1(x^1,x^2), y^2=f^2(x^1,x^2)$ satisfying the first order involutive system of finite Lie equations  $ y^1_2=0, \hspace{3mm} y^2_2-y^1_1=0 $ with no CC. As we have an algebraic Lie pseudogroup, a tricky computation ([24], p 467) allows to prove that $\Gamma $ is made by the transformations preserving a mixed tensor with square equal to zero as follows:\\
\[   \left( \begin{array}{rl}
y^1_1  & y^1_2 \\
y^2_1  & y^2_2
\end{array} 
\right) ^{-1}   
\left( \begin{array}{rl}
0  &  0  \\
1  &  0
\end{array}
\right)
 \left( \begin{array}{rl}
y^1_1  & y^1_2 \\
y^2_1  & y^2_2
\end{array} 
\right)  
=
\left( \begin{array}{rl}
0  &  0  \\
1  &  0
\end{array}
\right)    \]
We get the Lie form ${\Phi}^1\equiv y^1_2/y^1_1=0, {\Phi}^2\equiv (y^1_1)^2/(y^1_1y^2_2-y^1_2y^2_1)=1$ and let the reader exhibit $\cal{F}$.\\
Finally, introducing similarly the abstract number $i$ such that $i^2=-1$, we get the Cauchy-Riemann system  $y^2_2-y^1_1=0, \hspace{3mm} y^1_2+y^2_1=0$ with no CC defining complex analytic transformations and the correponding geometric object or {\it complex structure} is a mixed tensor with square equal to minus the $2\times 2$ identity matrix as we have now:     \\
\[   \left( \begin{array}{rl}
y^1_1  & y^1_2 \\
y^2_1  & y^2_2
\end{array} 
\right) ^{-1}   
\left( \begin{array}{rl}
0  &  -1  \\
1  &  0
\end{array}
\right)
 \left( \begin{array}{rl}
y^1_1  & y^1_2 \\
y^2_1  & y^2_2
\end{array} 
\right)  
=
\left( \begin{array}{rl}
0  &  -1  \\
1  &  0
\end{array}
\right)    \]

\noindent
{\bf EXAMPLE} 3.10: ({\it Riemann structure}) If $\omega$ is a section of ${\cal{F}}=S_2T^*$ with $det(\omega)\neq 0$ we get:\\
\noindent
{\it Lie form}  \hspace{25mm} ${\Phi}_{ij}(y_1)\equiv {\omega}_{kl}(y)y^k_iy^l_j={\omega}_{ij}(x)$\\
\noindent
{\it Medolaghi form} $  \hspace{2cm} {\Omega}_{ij}\equiv ({\cal{L}}(\xi)\omega)_{ij}\equiv {\omega}_{rj}(x){\partial}_i{\xi}^r+{\omega}_{ir}(x){\partial}_j{\xi}^r+{\xi}^r{\partial}_r{\omega}_{ij}(x)=0 $ \\
also called {\it Killing system} for historical reasons. A special section could be the Euclidean metric when $n=1,2,3$ as in elasticity theory or the Minkowski metric when $n=4$ as in special relativity. The main problem is that this system is not involutive unless we prolong the system to order two by differentiating once the equations. For such a purpose, introducing 
${\omega}^{-1}=({\omega}^{ij})$ as usual, we may define:\\
\noindent
{\it Christoffel symbols}$  \hspace{2cm}   {\gamma}^k_{ij}(x)=\frac{1}{2}{\omega}^{kr}(x)({\partial}_i{\omega}_{rj}(x) +{\partial}_j  {\omega}_{ri}(x) -{\partial}_r{\omega}_{ij}(x))={\gamma}^k_{ji}(x)  $\\
This is a new geometric object of order $2$ allowing to obtain, as in Example 3.5, an isomorphism $j_1(\omega)\simeq (\omega,\gamma)$ and the second order equations with $f^{-1}_1=g_1$:\\
\noindent
{\it Lie form } $ \hspace{2cm}  g^k_l(y^l_{ij}+{\gamma}^l_{rs}(y)y^r_iy^s_j)={\gamma}^k_{ij}(x)  $\\
\noindent
{\it Medolaghi form} $ \hspace{2mm} {\Gamma}^k_{ij}\equiv ({\cal{L}}(\xi)\gamma)^k_{ij}\equiv {\partial}_{ij}{\xi}^k+{\gamma}^k_{rj}(x){\partial}_i{\xi}^r+{\gamma}^k_{ir}(x){\partial}_j{\xi}^r-{\gamma}^r_{ij}(x){\partial}_r{\xi}^k+{\xi}^r{\partial}_r{\gamma}^k_{ij}(x)=0$\\
where $({\Gamma}^k_{ij})$ is a section of $S_2T^*\otimes T$. Surprisingly, the following expression:\\
\noindent
{\it Riemann tensor} $ \hspace{2cm}{\rho}^k_{lij}(x)\equiv {\partial}_i{\gamma}^k_{lj}(x)-{\partial}_j{\gamma}^k_{li}(x)+{\gamma}^r_{lj}(x){\gamma}^k_{ri}(x)-{\gamma}^r_{li}(x){\gamma}^k_{rj}(x)  $\\
is still a first order geometric object and even a tensor as a section of ${\wedge}^2T^*\otimes T^*\otimes T$ satisfying the purely algebraic relations :\\
    \hspace*{2cm}${\rho}^k_{lij}+{\rho}^k_{ijl}+{\rho}^k_{jli}=0, \hspace{5mm} {\omega}_{rl}{\rho}^l_{kij}+{\omega}_{kr}{\rho}^r_{lij}=0  \Rightarrow {\rho}_{klij}={\omega}_{kr}{\rho}^r_{lij}={\rho}_{ijkl}$. \\
Accordingly, the IC must express that the new first order equations $({\cal{L}}(\xi)\rho)^k_{lij}=0$ are only linear combinations of the previous ones and we get the Vessiot structure equations:\\
 \hspace*{4cm}${\rho}^k_{lij}(x)=c({\delta}^k_i{\omega}_{lj}(x)-{\delta}^k_j{\omega}_{li}(x))$\\
 describing the constant Riemannian curvature condition of Eisenhart [10]. Finally, as we have ${\rho}^r_{rij}(x)={\partial}_i{\gamma}^r_{rj}(x)-{\partial}_j{\gamma}^r_{ri}(x)=0$, we can only introduce the {\it Ricci tensor} ${\rho}_{ij}(x)={\rho}^r_{irj}(x)={\rho}_{ji}(x)$ by contracting indices and the scalar curvature $\rho(x)={\omega}^{ij}(x){\rho}_{ij}(x)$ in order to obtain $\rho(x)=n(n-1)c$. It remains to obtain all these results in a purely formal way, for example to prove that the number of components of the Riemann tensor is equal to $n^2(n^2-1)/12$ without dealing with indices.\\
 
 \noindent
{\bf REMARK} 3.11: Comparing the various Vessiot structure equations containing structure constants, we discover at once that the many $c$ appearing in the MC equations are {\it absolutely on equal footing} with the only $c$ appearing in the other examples. As their factors are either constant, linear or quadratic, {\it any identification of the quadratic terms appearing in the Riemann tensor with the quadratic terms appearing in the MC equations is definitively not correct} or, in an equivalent but more abrupt way, {\it the Cartan structure equations have nothing to do with the Vessiot structure equations}. As we shall see, {\it most of mathematical physics today is based on such a confusion}.\\

\noindent
{\bf REMARK} 3.12: Let us consider again Example 3.5 with ${\partial}_{xx}f(x)/{\partial}_xf(x)=\bar{\omega}(x)$ and introduce a variation $\eta(f(x))=\delta f(x)$ as in analytical or continuum mechanics. We get similarly $\delta {\partial}_xf={\partial}_x\delta f=\frac{\partial \eta}{\partial y}{\partial}_xf$ and so on, a result leading to $\delta \bar{\omega}(x)={\partial}_xf {\cal{L}}(\eta)\omega (f(x))$ where the Lie derivative involved is computed {\it over the target}. Let us now {\it pass from the target to the source} by introducing $\eta=\xi{\partial}_xf \Rightarrow \frac{\partial \eta}{\partial y}{\partial}_xf={\partial}_x\xi{\partial}_xf+\xi{\partial}_{xx}f$ and so on, a result leading to the particularly simple variation $\delta \bar{\omega}={\cal{L}}(\xi)\bar{\omega}$ {\it over the soure}. As another example of this general variational procedure, let us compare with the similar variations on which classical finite elasticity theory is based. Starting now with 
${\omega}_{kl}(f(x)){\partial}_if^k(x){\partial}_jf^l(x)={\bar{\omega}}_{ij}(x)$, where $\omega$ is the Euclidean metric, we obtain $(\delta{\bar{\omega}})_{ij}(x)={\partial}_if^k(x){\partial}_jf^l(x)({\cal{L}}(\eta)\omega)_{kl}(f(x))$ where the Lie derivative involved is computed {\it over the target}. {\it Passing now from the target to the source} as before, we find the particularly simple variation $\delta \bar{\omega}={\cal{L}}(\xi)\bar{\omega}$ {\it over the source}. For "small" deformations, source and target are of course identified but it is not true that the infinitesimal deformation tensor is in general the limit of the finite deformation tensor (for a counterexample, see [25], p 70). \\
  
Introducing a copy $Y$ of $X$ in the general framework, $(f,\delta f) $ must be considered as a section of 
$V(X\times Y)=(X\times Y){\times}_YT(Y)=X\times T(Y)$ over $X$. When $f$ is invertible (care), then we may consider the map $f:X\rightarrow Y:(x)\rightarrow (y=f(x))$ and define $\xi\in T$ by $\eta=T(f)(\xi)$ or rather $\eta=j_1(f)(\xi)$ in the language of geometric object, as a way to identify $f^{-1}(V(X\times Y))$ with $T=T(X)$. When $f=id$, this identification is canonical by considering vertical vectors along the diagonal 
$\Delta=\{(x,y)\in X\times Y{\mid} y=x\}$ and we get $\delta\omega=\Omega\in F_0={\omega}^{-1}(V({\cal{F}}))$. We point out that the above {\it vertical procedure} is a nice tool for studying nonlinear systems ([26], III, C and [27], III, 2).\\

\noindent
{\bf 4   JANET VERSUS SPENCER : THE LINEAR SEQUENCES}\\

 Let $\mu=({\mu}_1,...,{\mu}_n)$ be a multi-index with {\it length} ${\mid}\mu{\mid}={\mu}_1+...+{\mu}_n$, {\it class} $i$ if ${\mu}_1=...={\mu}_{i-1}=0,{\mu}_i\neq 0$ and $\mu +1_i=({\mu}_1,...,{\mu}_{i-1},{\mu}_i +1, {\mu}_{i+1},...,{\mu}_n)$. We set $y_q=\{y^k_{\mu}{\mid} 1\leq k\leq m, 0\leq {\mid}\mu{\mid}\leq q\}$ with $y^k_{\mu}=y^k$ when ${\mid}\mu{\mid}=0$. If $E$ is a vector bundle over $X$ with local coordinates $(x^i,y^k)$ for $i=1,...,n$ and $k=1,...,m$, we denote by $J_q(E)$ the $q$-{\it jet bundle} of $E$ with local coordinates simply denoted by $(x,y_q)$ and {\it sections} $f_q:(x)\rightarrow (x,f^k(x), f^k_i(x), f^k_{ij}(x),...)$ transforming like the section $j_q(f):(x)\rightarrow (x,f^k(x),{\partial}_if^k(x),{\partial}_{ij}f^k(x),...)$ when $f$ is an arbitrary section of $E$. Then both $f_q\in J_q(E)$ and $j_q(f)\in J_q(E)$ are over $f\in E$ and the {\it Spencer operator} just allows to distinguish them by introducing a kind of "{\it difference}" through the operator $D:J_{q+1}(E)\rightarrow T^*\otimes J_q(E): f_{q+1}\rightarrow j_1(f_q)-f_{q+1}$ with local components $({\partial}_if^k(x)-f^k_i(x), {\partial}_if^k_j(x)-f^k_{ij}(x),...) $ and more generally $(Df_{q+1})^k_{\mu,i}(x)={\partial}_if^k_{\mu}(x)-f^k_{\mu+1_i}(x)$. In a symbolic way, {\it when changes of coordinates are not involved}, it is sometimes useful to write down the components of $D$ in the form $d_i={\partial}_i-{\delta}_i$ and the restriction of $D$ to the kernel $S_{q+1}T^*\otimes E$ of the canonical projection ${\pi}^{q+1}_q:J_{q+1}(E)\rightarrow J_q(E)$ is {\it minus} the {\it Spencer map} $\delta=dx^i\wedge {\delta}_i:S_{q+1}T^*\otimes E\rightarrow T^*\otimes S_qT^*\otimes E$. The kernel of $D$ is made by sections such that $f_{q+1}=j_1(f_q)=j_2(f_{q-1})=...=j_{q+1}(f)$. Finally, if $R_q\subset J_q(E)$ is a {\it system} of order $q$ on $E$ locally defined by linear equations ${\Phi}^{\tau}(x,y_q)\equiv a^{\tau\mu}_k(x)y^k_{\mu}=0$ and local coordinates $(x,z)$ for the parametric jets up to order $q$, the $r$-{\it prolongation} $R_{q+r}={\rho}_r(R_q)=J_r(R_q)\cap J_{q+r}(E)\subset J_r(J_q(E))$ is locally defined when $r=1$ by the linear equations ${\Phi}^{\tau}(x,y_q)=0, d_i{\Phi}^{\tau}(x,y_{q+1})\equiv a^{\tau\mu}_k(x)y^k_{\mu+1_i}+{\partial}_ia^{\tau\mu}_k(x)y^k_{\mu}=0$ and has {\it symbol} $g_{q+r}=R_{q+r}\cap S_{q+r}T^*\otimes E\subset J_{q+r}(E)$ if one looks at the {\it top order terms}. If $f_{q+1}\in R_{q+1}$ is over $f_q\in R_q$, differentiating the identity $a^{\tau\mu}_k(x)f^k_{\mu}(x)\equiv 0$ with respect to $x^i$ and substracting the identity $a^{\tau\mu}_k(x)f^k_{\mu+1_i}(x)+{\partial}_ia^{\tau\mu}_k(x)f^k_{\mu}(x)\equiv 0$, we obtain the identity $a^{\tau\mu}_k(x)({\partial}_if^k_{\mu}(x)-f^k_{\mu+1_i}(x))\equiv 0$ and thus the restriction $D:R_{q+1}\rightarrow T^*\otimes R_q$ ([23],[27],[33]). \\
    
\noindent
{\bf DEFINITION} 4.1: $R_q$ is said to be {\it formally integrable} when the restriction ${\pi}^{q+1}_q:R_{q+1}\rightarrow R_q $ is an epimorphism $\forall r\geq 0$ or, equivalently, when all the equations of order 
$q+r$ are obtained by $r$ prolongations only $\forall r\geq 0$. In that case, $R_{q+1}\subset J_1(R_q)$ is a canonical equivalent formally integrable first order system on $R_q$ with no zero order equations, called the {\it Spencer form}.\\

\noindent
{\bf DEFINITION} 4.2: $R_q$ is said to be {\it involutive} when it is formally integrable and all the sequences $... \stackrel{\delta}{\rightarrow} {\wedge}^sT^*\otimes g_{q+r}\stackrel{\delta}{\rightarrow}...$ are exact $\forall 0\leq s\leq n, \forall r\geq 0$. Equivalently, using a linear change of local coordinates if necessary, we may {\it successively} solve the maximum number ${\beta}^n_q, {\beta}^{n-1}_q, ... , {\beta}^1_q$ of equations with respect to the principal jet coordinates of strict order $q$ and class $n,n-1,...,1$ in order to introduce the {\it characters} ${\alpha}^i_q=m\frac{(q+n-i-1)!}{(q-1)!((n-i)!}-{\beta}^i_q$ for $i=1, ..., n$ with ${\alpha}^n_q=\alpha$. Then $R_q$ is involutive if $R_{q+1}$ is obtained by only prolonging the ${\beta}^i_q$ equations of class $i$ with respect to $d_1,...,d_i$ for $i=1,...,n$. In that case $dim(g_{q+1})={\alpha}^1_q+...+{\alpha}^n_q$ and one can exhibit the {\it Hilbert polynomial} $dim(R_{q+r})$ in $r$ with leading term $(\alpha/n!)r^n$ when $\alpha \neq 0$. Such a prolongation procedure allows to compute {\it in a unique way} the principal ($pri$) jets from the parametric ($par$) other ones. This definition may also be applied to nonlinear systems as well.\\

  We obtain the following theorem generalizing for PD control systems the well known first order Kalman form of OD control systems where the derivatives of the input do not appear ([27], VI,1.14, p 802):\\

\noindent
{\bf THEOREM} 4.3: When $R_q$ is involutive, its Spencer form is involutive and can be modified to a {\it reduced Spencer form} in such a way that $\beta=dim(R_q)-\alpha$ equations can be solved with respect to the jet coordinates $z^1_n,...,z^{\beta}_n$ while $z^{\beta+1}_n,...,z^{\beta+\alpha}_n$ do not appear. In this case $z^{\beta+1},...,z^{\beta+\alpha}$ do not appear in the other equations.\\

   When $R_q$ is involutive, the linear differential operator ${\cal{D}}:E\stackrel{j_q}{\rightarrow} J_q(E)\stackrel{\Phi}{\rightarrow} J_q(E)/R_q=F_0$ of order $q$ with space of solutions $\Theta\subset E$ is said to be {\it involutive} and one has the canonical {\it linear Janet sequence} ([4], p 144):\\
\[  0 \longrightarrow  \Theta \longrightarrow T \stackrel{\cal{D}}{\longrightarrow} F_0 \stackrel{{\cal{D}}_1}{\longrightarrow}F_1 \stackrel{{\cal{D}}_2}{\longrightarrow} ... \stackrel{{\cal{D}}_n}{\longrightarrow} F_n \longrightarrow 0   \]
where each other operator is first order involutive and generates the {\it compatibility conditions} (CC) of the preceding one. As the Janet sequence can be cut at any place, {\it the numbering of the Janet bundles has nothing to do with that of the Poincar\'{e} sequence}, contrary to what many physicists  believe.\\

\noindent
{\bf DEFINITION} 4.4: The Janet sequence is said to be {\it locally exact at} $F_r$ if any local section of $F_r$ killed by ${\cal{D}}_{r+1}$ is the image by ${\cal{D}}_r$ of a local section of $F_{r-1}$. It is called {\it locally exact} if it is locally exact at each $F_r$ for $0\leq r \leq n$. The Poincar\'{e} sequence is locally exact but counterexemples may exist ([23], p 202).\\

    Equivalently, we have the involutive {\it first Spencer operator} $D_1:C_0=R_q\stackrel{j_1}{\rightarrow}J_1(R_q)\rightarrow J_1(R_q)/R_{q+1}\simeq T^*\otimes R_q/\delta (g_{q+1})=C_1$ of order one induced by $D:R_{q+1}\rightarrow T^*\otimes R_q$. Introducing the {\it Spencer bundles} $C_r={\wedge}^rT^*\otimes R_q/{\delta}({\wedge}^{r-1}T^*\otimes g_{q+1})$, the first order involutive ($r+1$)-{\it Spencer operator} $D_{r+1}:C_r\rightarrow C_{r+1}$ is induced by $D:{\wedge}^rT^*\otimes R_{q+1}\rightarrow {\wedge}^{r+1}T^*\otimes R_q:\alpha\otimes {\xi}_{q+1}\rightarrow d\alpha\otimes {\xi}_q+(-1)^r\alpha\wedge D{\xi}_{q+1}$ and we obtain the canonical {\it linear Spencer sequence} ([4], p 150):\\
\[    0 \longrightarrow \Theta \stackrel{j_q}{\longrightarrow} C_0 \stackrel{D_1}{\longrightarrow} C_1 \stackrel{D_2}{\longrightarrow} C_2 \stackrel{D_3}{\longrightarrow} ... \stackrel{D_n}{\longrightarrow} C_n\longrightarrow 0  \]
\noindent
as the Janet sequence for the first order involutive system $R_{q+1}\subset J_1(R_q)$.\\
    The Janet sequence and the Spencer sequence are connected by the following {\it crucial} commutative {\it diagram} (1) where the Spencer sequence is induced by the locally exact central horizontal sequence which is at the same time the Janet sequence for $j_q$ and the Spencer sequence for $J_{q+1}(E)\subset J_1(J_q(E))$ ([25], p 152):\\
  \[  SPENCER \hspace{5mm} SEQUENCE   \]  
 \[  \begin{array}{rcccccccccccl}
 &&&&& 0 &&0&&0&  &0&  \\
 &&&&& \downarrow && \downarrow && \downarrow &    & \downarrow &  \\
  & 0& \longrightarrow& \Theta &\stackrel{j_q}{\longrightarrow}&C_0 &\stackrel{D_1}{\longrightarrow}& C_1 &\stackrel{D_2}{\longrightarrow} & C_2 &\stackrel{D_3}{\longrightarrow} ... \stackrel{D_n}{\longrightarrow}& C_n &\longrightarrow 0 \\
  &&&&& \downarrow & & \downarrow & & \downarrow & &\downarrow &     \\
   & 0 & \longrightarrow & E & \stackrel{j_q}{\longrightarrow} & C_0(E) & \stackrel{D_1}{\longrightarrow} & C_1(E) &\stackrel{D_2}{\longrightarrow} & C_2(E) &\stackrel{D_3}{\longrightarrow} ... \stackrel{D_n}{\longrightarrow} & C_n(E) &   \longrightarrow 0 \\
   & & & \parallel && \hspace{5mm}\downarrow {\Phi}_0 & &\hspace{5mm} \downarrow {\Phi}_1 & & \hspace{5mm}\downarrow {\Phi}_2 &  & \hspace{5mm}\downarrow {\Phi}_n & \\
   0 \longrightarrow & \Theta &\longrightarrow & E & \stackrel{\cal{D}}{\longrightarrow} & F_0 & \stackrel{{\cal{D}}_1}{\longrightarrow} & F_1 & \stackrel{{\cal{D}}_2}{\longrightarrow} & F_2 & \stackrel{{\cal{D}}_3}{\longrightarrow} ... \stackrel{{\cal{D}}_n}{\longrightarrow} & F_n & \longrightarrow  0 \\
   &&&&& \downarrow & & \downarrow & & \downarrow &   &\downarrow &   \\
   &&&&& 0 && 0 && 0 &&0 &  
   \end{array}     \]
   \[    \hspace{5mm}  JANET  \hspace{5mm} SEQUENCE   \]
\vspace{5mm}  \\
\noindent   
In this diagram, {\it only depending on the left commutative square} ${\cal{D}}=\Phi\circ j_q$, the epimorhisms ${\Phi}_r:C_r(E)\rightarrow F_r$ for $0\leq r \leq n$ are successively induced by the canonical projection $\Phi={\Phi}_0:C_0(E)=J_q(E)\rightarrow J_q(E)/R_q=F_0$.\\
 
\noindent
{\bf EXAMPLE} 4.5 : ({\it Screw structure}): The system $R_1\subset J_1(T)$ defined by ${\xi}^1_2=0, {\xi}^2_2-{\xi}^1_1=0$ is involutive with $par(R_2)=\{{\xi}^1,{\xi}^2, {\xi}^1_1,{\xi}^2_1,{\xi}^1_{11},{\xi}^2_{11}\}$. The Spencer operator is not involutive as it is not even formally integrable because ${\partial}_2{\xi}^2_1-{\xi}^1_{11}=0, {\partial}_1{\xi}^2_1-{\xi}^2_{11}=0\Rightarrow {\partial}_1{\xi}^1_{11}-{\partial}_2{\xi}^2_{11}=0$. We obtain $dim(F_0)=2,dim(C_0(T))=6\Rightarrow dim(C_0)=dim(R_1)=4, dim(F_1)=0\Rightarrow dim(C_1(T))=dim(C_1)=6, dim(C_2(T))=dim(C_2)=2$ and it is not evident at all that the first order involutive operator $D_1:C_0\rightarrow C_1$ is defined by the $6$ PD equations:\\
\[   {\partial}_2{\xi}^1=0,{\partial}_2{\xi}^2-{\xi}^1_1=0, {\partial}_2{\xi}^1_1=0, {\partial}_2{\xi}^2_1-{\partial}_1{\xi}^1_1=0, {\partial}_1{\xi}^1-{\xi}^1_1=0, {\partial}_1{\xi}^2-{\xi}^2_1=0\]
The case of a complex structure is similar and left to the reader.\\

\noindent
{\bf 5   DIFFERENTIAL MODULES AND INVERSE SYSTEMS}\\

   An important but difficult problem in engineering physics is to study how the formal properties of a system of order $q$ with $n$ independent variables and $m$ unknowns depend on the parameters involved in that system. This is particularly clear in classical control theory where the systems are classified into two categories, namely the "controllable" ones and the "uncontrollable" ones ([14],[27]). In order to understand the problem studied by Macaulay in [M], that is roughly to determine the minimum number of  solutions of a system that must be known in order to determine all the others by using derivatives and linear combinations with constant coefficients in a field $k$, let us start with the following motivating example:\\

\noindent
{\bf EXAMPLE} 5.1:  When $n=1,m=1, q=3$, using a sub-index $x$ for the derivatives with $d_xy=y_x$ and so on, the general solution of $y_{xxx}-y_x=0$ is $y=a e^x+ b e^{-x}+c 1$ with $a,b,c$ constants and the derivative of $e^x$ is $e^x$, the derivative of $e^{-x}$ is $-e^{-x}$ and the derivative of $1$ is $0$. Hence we could believe that we need a basis $\{1,e^x, e^{-x}\}$ with {\it three} generators for obtaining all the solutions through derivatives. Also, when $n=1, m=2, k=\mathbb{R}$ and $a$ is a constant real parameter, the OD system $y^1_{xx}-ay^1=0, y^2_x=0$ needs two generators $\{(x,0),(0,1)\}$ when $a=0$ with the only $d_x$ killing both $y^1_x$ and $y_2$ but only one generator when $a\neq 0$, namely $\{(ch(x),1)\} $ when $a=1$. Indeed, setting $y=y^1-y^2$ brings $y^1=y_{xx}, y^2=y_{xx}-y$ and an equivalent system defined by the single OD equation $y_{xxx}-y_x=0$ for the only $y$. Introducing the corresponding poynomial ideal $({\chi}^3-\chi)=(\chi)\cap ({\chi}-1)\cap (\chi +1)$, we check that $d_x$ kills $y_{xx}-y$, $d_x-1$ kills $y_{xx}+y_x$ and $d_x+1$ kills $y_{xx}-y_x$, a result leading, as we shall see, to the {\it only} generator $\{ch(x)-1\}$.\\

   More precisely, if $K$ is a differential field containing $\mathbb{Q}$ with $n$ commuting {\it derivations} 
   ${\partial}_i$, that is to say ${\partial}_i(a+b)={\partial}_ia+{\partial}_ib$ and ${\partial}_i(ab)=({\partial}_ia)b+a{\partial}_ib, \forall a,b\in K$ for $i=1,...,n$, we denote by $k$ a subfield of constants. Let us introduce $m$ {\it differential indeterminates} $y^k$ for $k=1,...,m$ and $n$ commuting {\it formal derivatives} $d_i$ with $d_iy^k_{\mu}=y^k_{\mu+1_i} $. We introduce the non-commutative {\it ring of differential operators} $D=K[d_1,...,d_n]=K[d]$ with $d_ia=ad_i+{\partial}_ia, \forall a\in K$ in the operator sense and the {\it differential module} $Dy=Dy^1+...+Dy^m$. If $\{{\Phi}^{\tau}=a^{\tau\mu}_ky^k_{\mu}\}$ is a finite number of elements in $Dy$ indexed by $\tau$, we may introduce the {\it differential module of equations} $I=D\Phi\subset Dy$ and the finitely generated {\it residual differential module} $M=Dy/I$.\\
   
   In the algebraic framework considered, {\it only two possible formal constructions can be obtained from} $M$ when $D=K[d]$, namely $hom_D(M,D)$ and $M^*=hom_K(M,K)$ ([3],[27],[32]).\\
 
 \noindent  
{\bf THEOREM} 5.2: $hom_D(M,D)$ is a {\it right} differential module that can be {\it converted} to a {\it left} differential module by introducing the {\it right} differential module structure of ${\wedge}^nT^*$. As a differential geometric counterpart, we get the {\it formal adjoint} of ${\cal{D}}$, namely $ad({\cal{D}}):{\wedge}^nT^*\otimes F^* \rightarrow {\wedge}^nT^*\otimes E^*$ usually constructed through an integration by parts and where $E^*$ is obtained from $E$ by inverting the local transition matrices, the simplest example being the way $T^*$ is obtained from $T$.\\

\noindent
{\bf REMARK} 5.3: Such a result explains why dual objects in physics and engineering are no longer tensors but tensor {\it densities}, with no reference to any variational calculus. For example the EM potential is a section of $T^*$ and the EM field is a section of ${\wedge}^2T^*$ while the EM induction is a section of ${\wedge}^4T^*\otimes {\wedge}^2T\simeq {\wedge}^2T^*$ and the EM current is a section of ${\wedge}^4T^*\otimes T\simeq {\wedge}^3T^*$. \\

   The filtration $D_0=K\subseteq D_1=K\oplus T \subseteq ... \subseteq D_q\subseteq ... \subseteq D$ of $D$ by the order of operators induces a filtration/inductive limit $0\subseteq M_0\subseteq M_1 \subseteq ... \subseteq M_q\subseteq ... \subseteq M$ and provides by duality {\it over} $K$ the projective limit $M^*=R\rightarrow ... \rightarrow R_q \rightarrow ... \rightarrow R_1\rightarrow R_0\rightarrow 0$ of formally integrable systems. As $D$ is generated by $K$ and $T=D_1/D_0$, we can define for any $f\in M^*$:\\
 \[ (af)(m)=af(m)=f(am), ({\xi}f)(m)={\xi}f(m)-f({\xi}m), \forall a\in K, \forall {\xi}=a^id_i\in T, \forall m\in M \]
 and check $d_ia=ad_i+{\partial}_ia, \xi\eta-\eta\xi=[\xi,\eta]$ in the operator sense by introducing the standard bracket of vector fields on $T$. Finally we get $(d_if)^k_{\mu}=(d_if)(y^k_{\mu})={\partial}_if^k_{\mu}-f^k_{\mu+1_i}$ in a coherent way.\\
 
 \noindent
 {\bf THEOREM} 5.4: $R=M^*$ has a structure of differential module induced by the Spencer operator.\\
 
 \noindent
 {\bf REMARK} 5.5: When $m=1$ and $D=k[d]$ is a commutative ring isomorphic to the polynomial ring $A=k[\chi]$ for the indeterminates ${\chi}_1,...,{\chi}_n$, this result {\it exactly} describes the {\it inverse system} of Macaulay with $-d_i={\delta}_i$ ([M], \S 59,60).\\
 
 \noindent
{\bf DEFINITION} 5.6: A {\it simple} module is a module having no other proper submodule than $0$. A {\it semi-simple} module is a direct sum of simple modules. When A is a commutative integral domain and $M$ a finitely generated module over $A$, the {\it socle} of $M$ is the largest semi-simple submodule of $M$, that is $soc(M)=\oplus soc_{\mathfrak{m}}(M)$ where $soc_{\mathfrak{m}}(M)$ is the direct sum of all the {\it isotypical} simple submodules of $M$ isomorphic to $A/{\mathfrak{m}}$ for $\mathfrak{m}\in max(A)$ the set of maximal proper ideals of $A$. The {\it radical} of a module is the intersection of all its maximum {\it proper} submodules. The quotient of a module by its radical is called the {\it top} and is a semi-simple module ([3]).\\

  The "{\it secret} " of Macaulay is expressed by the next theorem:\\
  
\noindent
{\bf THEOREM} 5.7: Instead of using the socle of $M$ over $A$, one may use duality over $k$ in order to deal with the short exact sequence $0 \rightarrow rad(R) \rightarrow R \rightarrow top(R) \rightarrow 0 $ where $top(R)$ is the dual of $soc(M)$.\\

   However, Nakayama's lemma ([3],[19],[32]) cannot be used in general unless $R$ is finitely generated over $k$ and thus over $D$. The main idea of Macaulay has been to overcome this difficulty by dealing only with {\it unmixed} ideals when $m=1$. As a generalization, one can state ([27]):\\
   
\noindent
{\bf DEFINITION} 5.8: One has the {\it purity filtration} $0=t_n(M)\subseteq ... \subseteq t_0(M)=t(M)\subseteq M$ where {\it any} involutive system of order $p$ defining $Dm$ is such that ${\alpha}^{n-r}_p=0, ... , {\alpha}^n_p=0$ when $m\in t_r(M)$ and $M$ is said to be $r$-{\it pure} if $t_r(M)=0, t_{r-1}(M)=M$. With $t(M)=\{m\in M\mid \exists 0\neq a\in A, am=0\}$ we say that $M$ is a $0$-pure or {\it torsion-free} module if $t(M)=0$ and a {\it torsion module} if $t(M)=M$.\\

\noindent
{\bf EXAMPLE} 5.9: With $n=2, q=2$, let us consider the involutive system $y_{(0,2)}\equiv y_{22}=0, y_{(1,1)}\equiv y_{12}=0$. Then $z'=y_1$ satisfies $z'_2=0$ while $z''=y_2$ satisfies $z''_2=0, z''_1=0$ and we have the filtration $0=t_2(M)\subset t_1(M)\subset t_0(M)=t(M)=M$ with $z'' \in t_1(M), z' \in t_0(M)$ but $z' \notin t_1(M)$. This classification of observables has never been applied to engineering systems like the ones to be found in magnetohydrodynamics (MHD) because the mathematics involved are not known.\\

\noindent
{\bf REMARK} 5.10: A standard result in commutative algebra allows to embed any torsion-free module into a free module ([32]). Such a property provides the possibility to {\it parametrize} the solution space of the corresponding system of OD/PD equations by a finite number of potential like arbitrary functions. For this, in order to test the possibility to parametrize a given operator ${\cal{D}}_1$, one may construct the adjoint operator $ad({\cal{D}}_1)$ and look for generating CC in the form of an operator $ad({\cal{D}})$. As $ad({\cal{D}})\circ ad({\cal{D}}_1)=ad({\cal{D}}_1\circ {\cal{D}})=0\Rightarrow {\cal{D}}_1\circ {\cal{D}}=0$, it only remains to check that the CC of $\cal{D}$ are generated by ${\cal{D}}_1$. When $n=1$ this result amounts to Kalman test and the fact that a classical OD control system is controllable if and only if it is parametrizable, a result showing that {\it controllability is an intrinsic structural property of a control system}, not depending on the choice of inputs and outputs contrary to a well established engineering tradition ([14],[27]). When $n=2$, the formal adjoint of the only CC for the deformation tensor has been used in the Introduction in order to parametrize the stress equation by means of the Airy function. This result is also valid for the non-commutative ring $D=K[d]$.  \\ 

\noindent
{\bf EXAMPLE} 5.11: With $K=\mathbb{Q} (x^1,x^2,x^3)$, infinitesimal contact transformations are defined by the system ${\partial}_2{\xi}^1-x^3{\partial}_2{\xi}^2+x^3{\partial}_1{\xi}^1-(x^3)^2{\partial}_1{\xi}^2-{\xi}^3=0, \hspace{3mm} {\partial}_3{\xi}^1-x^3{\partial}_3{\xi}^2=0$. Multiplying by test functions $({\lambda}^1,{\lambda}^2)$ and integrating by parts, we obtain the adjoint operator (up to sign): \\
\[  {\partial}_2{\lambda}^1+x^3{\partial}_1{\lambda}^1+{\partial}_3{\lambda}^2={\mu}^1, \hspace{3mm}-x^3{\partial}_2{\lambda}^1-(x^3)^2{\partial}_1{\lambda}^1-x^3{\partial}_3{\lambda}^2-{\lambda}^2={\mu}^2, \hspace{3mm}{\lambda}^1={\mu}^3  \]
It follows that ${\lambda}^1={\mu}^3, {\lambda}^2=-{\mu}^2-x^3{\mu}^1\Rightarrow {\partial}_2{\mu}^3+x^3{\partial}_1{\mu}^3-{\partial}_3{\mu}^2-x^3{\partial}_3{\mu}^1-2{\mu}^1=0$. Multiplying again by a test function $\phi$, we discover the parametrization ${\xi}^1=x^3{\partial}_3\phi-\phi, {\xi}^2={\partial}_3\phi, {\xi}^3=-{\partial}_2\phi-x^3{\partial}_1\phi$ which is not evident at first sight.  \\

   When $M$ is $r$-pure, Theorem 4.3 provides the exact sequence $0\rightarrow M\rightarrow k({\chi}_1,...,{\chi}_{n-r})\otimes M$, also discovered by Macaulay ([M], \S 77, 82), and one obtains the following key result for studying the {\it identifiability} of OD/PD control systems (see {\it localization} in ([19],[27],32[29],[30],[32]).\\

\noindent
{\bf THEOREM} 5.12: When $M$ is $n$-pure, one may use the {\it chinese remainder theorem} ([19], p 41) in order to prove that the minimum number of generators of $R$ is equal to the maximum number of isotypical components that can be found among the various components of $soc(M)$ or $top(R)$.\\

\noindent
{    }  \\

\noindent
{\bf 6   JANET VERSUS SPENCER : THE NONLINEAR SEQUENCES}\\

Nonlinear operators do not in general admit CC as can be seen by considering the involutive example 
$y_{22}-\frac{1}{3}(y_{11})^3=u, y_{12}-\frac{1}{2}(y_{11})^2=v$ with $m=1,n=2,q=2$, contrary to what happens in the study of Lie pseudogroups. However, the kernel of a linear operator ${\cal{D}}:E\rightarrow F$ is always taken with respet to the zero section of $F$, while it must be taken with respect to a prescribed section by a {\it double arrow} for a nonlinear operator. Keeping in mind the linear Janet sequence and the examples of Vessiot structure equations already presented, one obtains:  \\

\noindent
{\bf THEOREM} 6.1: There exists a {\it nonlinear Janet sequence} associated with the Lie form of an involutive system of finite Lie equations:   \\
\[  \begin{array}{rcccl}
  & {\Phi}_{\omega}\circ j_q &   &  I\circ j_1  &   \\
  0\rightarrow \Gamma \rightarrow aut(X) &\rightrightarrows &  {\cal{F}}  &\rightrightarrows   &  {\cal{F}}_1 \\
     & \omega\circ\alpha  &  &  0 &
  \end{array}  \]
where the kernel of the first operator $f\rightarrow {\Phi}_{\omega}\circ j_q(f)={\Phi}_{\omega}(j_q(f))=j_q(f)^{-1}(\omega)$ is taken with respect to the section $\omega$ of $\cal{F}$ while the kernel of the second operator $\omega\rightarrow I(j_1(\omega))\equiv A(\omega){\partial}_x\omega+B(\omega)$ is taken with respect to the zero section of the vector bundle ${\cal{F}}_1$ over ${\cal{F}}$.\\

\noindent
{\bf COROLLARY} 6.2: By linearization at the identity, one obtains the involutive {\it Lie operator} ${\cal{D}}={\cal{D}}_{\omega}:T\rightarrow F_0:\xi\rightarrow {\cal{L}}(\xi)\omega $ with kernel $\Theta=\{\xi\in T{\mid}{\cal{L}}(\xi)\omega=0\}\subset T$ satisfying $[\Theta,\Theta]\subset \Theta$ and the corresponding {\it linear Janet sequence} where $F_0={\omega}^{-1}(V({\cal{F}}))$ and $F_1={\omega}^{-1}({\cal{F}}_1)$.\\

Now we notice that $T$ is a natural vector bundle of order $1$ and $J_q(T)$ is thus a natural vector bundle of order $q+1$. Looking at the way a vector field and its derivatives are transformed under any $f\in aut(X)$ while replacing $j_q(f)$ by $f_q$, we obtain:\\
\[  {\eta}^k(f(x))=f^k_r(x){\xi}^r(x) \Rightarrow {\eta}^k_u(f(x))f^u_i(x)=f^k_r(x){\xi}^r_i(x)+f^k_{ri}(x){\xi}^r(x)\]
and so on, a result leading to:\\

\noindent
{\bf LEMMA} 6.3: $J_q(T)$ is {\it associated} with ${\Pi}_{q+1}={\Pi}_{q+1}(X,X)$ that is we can obtain a new section ${\eta}_q=f_{q+1}({\xi}_q)$ from any section ${\xi}_q \in J_q(T)$ and any section $f_{q+1}\in {\Pi}_{q+1}$ by the formula:\\
\[ d_{\mu}{\eta}^k\equiv {\eta}^k_rf^r_{\mu}+ ...=f^k_r{\xi}^r_{\mu}+  ...  +f^k_{\mu+1_r}{\xi}^r , \forall 0\leq {\mid}\mu {\mid}\leq q\]
where the left member belongs to $V({\Pi}_q)$. Similarly $R_q\subset J_q(T)$ is associated with ${\cal{R}}_{q+1}\subset {\Pi}_{q+1}$.\\

In order to construct another nonlinear sequence, we need a few basic definitions on {\it Lie groupoids} and {\it Lie algebroids} that will become substitutes for Lie groups and Lie algebras. As in the beginning of section 3, the first idea is to use the chain rule for derivatives $j_q(g\circ f)=j_q(g)\circ j_q(f)$ whenever $f,g\in aut(X)$ can be composed and to replace both $j_q(f)$ and $j_q(g)$ respectively by $f_q$ and $g_q$ in order to obtain the new section $g_q\circ f_q$. This kind of "composition" law can be written in a pointwise symbolic way by introducing another copy $Z$ of $X$ with local coordinates $(z)$ as follows:\\
\[ {\gamma}_q:{\Pi}_q(Y,Z){\times}_Y{\Pi}_q(X,Y)\rightarrow {\Pi}_q(X,Z):((y,z,\frac{\partial z}{\partial y},...),(x,y,\frac{\partial y}{\partial x},...)\rightarrow (x,z,\frac{\partial z}{\partial y}\frac{\partial y}{\partial x},...)      \]
We may also define $j_q(f)^{-1}=j_q(f^{-1})$ and obtain similarly an "inversion" law.\\

\noindent
{\bf DEFINITION} 6.4: A fibered submanifold ${\cal{R}}_q\subset {\Pi}_q$ is called a {\it system of finite Lie equations} or a {\it Lie groupoid} of order $q$ if we have an induced {\it source projection} ${\alpha}_q:{\cal{R}}_q\rightarrow X$, {\it target projection} ${\beta}_q:{\cal{R}}_q\rightarrow X$, {\it composition} ${\gamma}_q:{\cal{R}}_q{\times}_X{\cal{R}}_q\rightarrow {\cal{R}}_q$, {\it inversion} ${\iota}_q:{\cal{R}}_q\rightarrow {\cal{R}}_q$ 
and {\it identity} $id_q:X\rightarrow {\cal{R}}_q$. In the sequel we shall only consider {\it transitive} Lie groupoids such that the map $({\alpha}_q,{\beta}_q):{\cal{R}}_q\rightarrow X\times X $ is an epimorphism and we shall denote by ${\cal{R}}^0_q=id^{-1}({\cal{R}}_q)$ the {\it isotropy Lie group bundle} of 
${\cal{R}}_q$. Also, one can prove that the new system ${\rho}_r({\cal{R}}_q)={\cal{R}}_{q+r}$ obtained by differentiating $r$ times all the defining equations of ${\cal{R}}_q$ is a Lie groupoid of order $q+r$. Finally, one can write down the Lie form and obtain ${\cal{R}}_q=\{ f_q\in{\Pi}_q{\mid} f^{-1}_q(\omega)=\omega\}$.\\

Now, using the {\it algebraic bracket} $\{ j_{q+1}(\xi),j_{q+1}(\eta)\}=j_q([\xi,\eta]), \forall \xi\eta\in T$, we may  obtain by bilinearity a {\it differential bracket} on $J_q(T)$ extending the bracket on $T$:\\
\[   [{\xi}_q,{\eta}_q]=\{{\xi}_{q+1},{\eta}_{q+1}\}+i(\xi)D{\eta}_{q+1}-i(\eta)D{\xi}_{q+1}, \forall {\xi}_q,{\eta}_q\in J_q(T) \]
which does not depend on the respective lifts ${\xi}_{q+1}$ and ${\eta}_{q+1}$ of ${\xi}_q$ and ${\eta}_q$ in $J_{q+1}(T)$. This bracket on sections satisfies the Jacobi identity and we set:\\

\noindent
{\bf DEFINITION} 6.5: We say that a vector subbundle $R_q\subset J_q(T)$ is a {\it system of infinitesimal Lie equations} or a {\it Lie algebroid} if $[R_q,R_q]\subset R_q$, that is to say $[{\xi}_q,{\eta}_q]\in R_q, \forall {\xi}_q,{\eta}_q\in R_q$. The kernel $R^0_q$ of the projection ${\pi}^q_0:R_q\rightarrow T$ is the {\it isotropy Lie algebra bundle} of ${\cal{R}}^0_q$ and $[R^0_q,R^0_q]\subset R^0_q$ does not contain derivatives.\\

\noindent
{\bf PROPOSITION} 6.6 : There is a nonlinear differential sequence:\\
\[ 0\longrightarrow aut(X) \stackrel{j_{q+1}}{\longrightarrow} {\Pi}_{q+1}(X,X)\stackrel{\bar{D}}{\longrightarrow}T^*\otimes J_q(T)\stackrel{{\bar{D}}'}{\longrightarrow} {\wedge}^2T^*\otimes J_{q-1}(T)  \]
with $\bar{D}f_{q+1}\equiv f_{q+1}^{-1}\circ j_1(f_q)-id_{q+1}={\chi}_q \Rightarrow {\bar{D}}'{\chi}_q(\xi,\eta)\equiv D{\chi}_q(\xi,\eta)-\{{\chi}_q(\xi),{\chi}_q(\eta)\}=0 $. Moreover, setting ${\chi}_0=A-id\in T^*\otimes T$, this sequence is locally exact if $det(A)\neq 0$.\\

\noindent
{\bf Proof}: There is a canonical inclusion ${\Pi}_{q+1}\subset J_1({\Pi}_q)$ defined by $y^k_{\mu,i}=y^k_{\mu+1_i}$ and the composition $f^{-1}_{q+1}\circ j_1(f_q)$ is a well defined section of $J_1({\Pi}_q)$ over the section $f^{-1}_q\circ f_q=id_q$ of ${\Pi}_q$ like $id_{q+1}$. The difference ${\chi}_q=f^{-1}_{q+1}\circ j_1(f_q)-id_{q+1}$ is thus a section of $T^*\otimes V({\Pi}_q)$ over $id_q$ and we have already noticed that 
$ id^{-1}_q(V({\Pi}_q))=J_q(T)$. For $q=1$ we get with $g_1=f^{-1}_1$:\\
\[ {\chi}^k_{,i}=g^k_l{\partial}_if^l-{\delta}^k_i=A^k_i-{\delta}^k_i,\hspace{5mm} {\chi}^k_{j,i}=g^k_l({\partial}_if^l_j-A^r_if^l_{rj})  \]
We also obtain from Lemma 6.3 the useful formula $ f^k_r{\chi}^r_{\mu,i}+...+f^k_{\mu+1_r}{\chi}^r_{,i}={\partial}_if^k_{\mu}-f^k_{\mu+1_i}$ allowing to determine ${\chi}_q$ inductively.\\
We refer to ([26], p 215) for the inductive proof of the local exactness, providing the only formulas that will be used later on and can be checked directly by the reader:\\
\[  {\partial}_i{\chi}^k_{,j}-{\partial}_j{\chi}^k_{,i}-{\chi}^k_{i,j}+{\chi}^k_{j,i}-({\chi}^r_{,i}{\chi}^k_{r,j}-{\chi}^r_{,j}{\chi}^k_{r,i})=0  \]
\[  {\partial}_i{\chi}^k_{l,j}-{\partial}_j{\chi}^k_{l,i}-{\chi}^k_{li,j}+{\chi}^k_{lj,i}-({\chi}^r_{,i}{\chi}^k_{lr,j}+{\chi}^r_{l,i}{\chi}^k_{r,j}-{\chi}^r_{l,j}{\chi}^k_{r,i}-{\chi}^r_{,j}{\chi}^k_{lr,i})=0\]
There is no need for double-arrows in this framework as the kernels are taken with respect to the zero section of the vector bundles involved. We finally notice that the main difference with the gauge sequence is that {\it all the indices range from} $1$ {\it to} $n$ and that the condition $det(A)\neq 0$ amounts to $\Delta=det({\partial}_if^k)\neq 0$ because $det(f^k_i)\neq 0$ by assumption. \hspace{1cm}
$ \Box $ \\

\noindent
{\bf COROLLARY} 6.7: There is a restricted nonlinear differential sequence:\\
\[ 0\longrightarrow \Gamma \stackrel{j_{q+1}}{\longrightarrow} {\cal{R}}_{q+1} \stackrel{\bar{D}}{\longrightarrow} T^*\otimes R_q\stackrel{{\bar{D}}'}{\longrightarrow}{\wedge}^2T^*\otimes J_{q-1}(T)  \]

\noindent
{\bf DEFINITION} 6.8: A {\it splitting} of the short exact sequence $0\rightarrow R^0_q\rightarrow R_q\stackrel{{\pi}^q_0}{\rightarrow} T \rightarrow 0$ is a map ${\chi}'_q:T\rightarrow R_q$ such that ${\pi}^q_0\circ {\chi}'_q=id_T$ or equivalently a section of $T^*\otimes R_q$ over $id_T\in T^*\otimes T$ and is called a $R_q$-{\it connection}. Its {\it curvature} ${\kappa}'_q\in {\wedge}^2T^*\otimes R^0_q$ is defined by ${\kappa}'_q(\xi,\eta)=[{\chi}'_q(\xi),{\chi}'_q(\eta)]-{\chi}'_q([\xi,\eta])$. We notice that ${\chi}'_q=-{\chi}_q$ is a connection with ${\bar{D}}'{\chi}'_q={\kappa}'_q$ if and only if $A=0$. In particular $({\delta}^k_i,-{\gamma}^k_{ij})$ is the only existing symmetric connection for the Killing system.         \\

\noindent
{\bf REMARK} 6.9: Rewriting the previous formulas with $A$ instead of ${\chi}_0$ we get:  \\
\[ {\partial}_iA^k_j-{\partial}_jA^k_i-A^r_i{\chi}^k_{r,j}+A^r_j{\chi}^k_{r,i}=0  \]
\[ {\partial}_i{\chi}^k_{l,j}-{\partial}_j{\chi}^k_{l,i}-{\chi}^r_{l,i}{\chi}^k_{r,j}+{\chi}^r_{l,j}{\chi}^k_{r,i}-A^r_i{\chi}^k_{lr,j}+A^r_j{\chi}^k_{lr,i}=0  \]
When $q=1, g_2=0$ and though surprising it may look like, we find back {\it exactly} all the formulas presented by E. and F. Cosserat in ([C], p 123 and [16]). Even more strikingly, {\it in the case of a Riemann structure, the last two terms disappear but the quadratic terms are left while, in the case of screw and complex structures, the quadratic terms disappear but the last two terms are left}.\\

\noindent
{\bf COROLLARY} 6.10: When $det(A)\neq 0$ there is a nonlinear {\it stabilized} sequence at order $q$:\\
\[ 0 \longrightarrow {aut(X)} \stackrel{j_q}{\longrightarrow} {\Pi}_q \stackrel{{\bar{D}}_1}{\longrightarrow} C_1(T) \stackrel{{\bar{D}}_2}{\longrightarrow} C_2(T)  \]
called {\it nonlinear Spencer sequence} where ${\bar{D}}_1$ and ${\bar{D}}_2$ are involutive and its restriction:ÊÊ\\
\[ 0\longrightarrow \Gamma \stackrel{j_q}{\longrightarrow} {\cal{R}}_q \stackrel{{\bar{D}}_1}{\longrightarrow} C_1 \stackrel{{\bar{D}}_2}{\longrightarrow} C_2   \]
is such that ${\bar{D}}_1$ and ${\bar{D}}_2$ are involutive whenever ${\cal{R}}_q$ is involutive.\\

\noindent
{\bf Proof}: With ${\mid}\mu{\mid}=q$ we have ${\chi}^k_{\mu,i}=-g^k_lA^r_if^l_{\mu+1_r}+terms(order \leq q)$. Setting ${\chi}^k_{\mu,i}=A^r_i{\tau}^k_{\mu,r}$, we obtain ${\tau}^k_{\mu,r}=-g^k_lf^l_{\mu+1_r}+terms(order\leq q)$ and ${\bar{D}}:{\Pi}_{q+1}\rightarrow T^*\otimes J_q(T)$ restrics to ${\bar{D}}_1:{\Pi}_q\rightarrow C_1(T)$. \\
Finally, setting $ A^{-1}=B=id-{\tau}_0$, we obtain successively:\\
\[  {\partial}_i{\chi}^k_{\mu,j}-{\partial}_j{\chi}^k_{\mu,i}+terms({\chi}_q) -(A^r_i{\chi}^k_{\mu+1_r,j}-A^r_j{\chi}^k_{\mu+1_r,i})=0ÊÊ\]
\[  B^i_rB^j_s({\partial}_i{\chi}^k_{\mu,j}-{\partial}_j{\chi}^k_{\mu,i})+terms ({\chi}_q)-({\tau}^k_{\mu+1_r,s}-{\tau}^k_{\mu+1_s,r})=0  \]
We obtain therefore $D{\tau}_{q+1}+terms({\tau}_q)=0$ and ${\bar{D}}':T^*\otimes J_q(T)\rightarrow {\wedge}^2T^*\otimes J_{q-1}(T)$ restricts to ${\bar{D}}_2:C_1(T)\rightarrow C_2(T)$. \\
In the case of Lie groups of transformations, the symbol of the involutive system $R_q$ {\it must} be $g_q=0$ providing an isomorphism ${\cal{R}}_{q+1}\simeq {\cal{R}}_q\Rightarrow R_{q+1}\simeq R_q$ and we have therefore $C_r={\wedge}^rT^*\otimes R_q$ for $ r=1,...,n$ in the linear Spencer sequence. \hspace{1cm}$\Box$  \\

\noindent
{\bf REMARK} 6.11: The passage from ${\chi}_q$ to ${\tau}_q$ is {\it exactly} the one done by E. and F. Cosserat in ([C], p 190). However, even if is a good idea to pass from the source to the target, the way they realize it is based on a {\it subtle misunderstanding} that we shall correct later on in Proposition 6.16.\\

If $f_{q+1},g_{q+1}\in {\Pi}_{q+1}$ and $f'_{q+1}=g_{q+1}\circ f_{q+1}$, we get:\\
\[ {\bar{D}}f'_{q+1}=f^{-1}_{q+1}\circ g^{-1}_{q+1}\circ j_1(g_q)\circ j_1(f_q)-id_{q+1}=f^{-1}_{q+1}\circ {\bar{D}}g_{q+1}\circ j_1(f_q)+{\bar{D}}f_{q+1} \]

\noindent
{\bf DEFINITION} 6.12: For any section $f_{q+1}\in {\cal{R}}_{q+1}$, the transformation:\\
\[   {\chi}_q \longrightarrow  {\chi}'_q=f^{-1}_{q+1}\circ {\chi}_q\circ j_1(f_q)+{\bar{D}}f_{q+1}  \]
is called a {\it gauge transformation} and exchanges the solutions of the {\it field equations} 
${\bar{D}}'{\chi}_q=0$.  \\

Introducing the {\it formal Lie derivative} on $J_q(T)$ by the formulas:\\
\[  L({\xi}_{q+1}){\eta}_q=\{{\xi}_{q+1},{\eta}_{q+1}\}+i(\xi)D{\eta}_{q+1}=[{\xi}_q,{\eta}_q]+i(\eta)D{\xi}_{q+1} \]
\[      (L(j_1({\xi}_{q+1})){\chi}_q)(\zeta)= L({\xi}_{q+1})({\chi}_q(\zeta))-{\chi}_q([\xi,\zeta])   \]
and passing to the limit with $f_{q+1}=id_{q+1}+t {\xi}_{q+1}+ ... $ for $t\rightarrow 0$ {\it over the source}, we get:\\

\noindent
{\bf LEMMA} 6.13: An {\it infinitesimal gauge transformation} has the form:  \\
\[        \delta {\chi}_q= D{\xi}_{q+1}+ L(j_1({\xi}_{q+1})){\chi}_q     \]

Passing again to the limit but now {\it over the target} with ${\chi}_q=\bar{D}f_{q+1}$ and $g_{q+1}=id_{q+1}+ t {\eta}_{q+1}+ ... $, we obtain the variation:  \\
\[         \delta {\chi}_q= f^{-1}_{q+1}\circ D{\eta}_{q+1}\circ j_1(f_q)     \]
 
\noindent
{\bf PROPOSITION} 6.14: The same variation is obtained whenever ${\eta}_{q+1}=f_{q+2}({\xi}_{q+1}+{\chi}_{q+1}(\xi))$ with ${\chi}_{q+1}=\bar{D}f_{q+2}$, a transformation which only depends on $j_1(f_{q+1})$ and is invertible if and only if $det(A)\neq 0$.\\
\noindent
{\bf Proof}: Choosing  $f_{q+1},g_{q+1},h_{q+1}\in {\cal{R}}_{q+1}$ such that $g_{q+1}\circ f_{q+1}=f_{q+1}\circ h_{q+1}$ and passing to the limits $g_{q+1}=id_{q+1}+t{\eta}_{q+1}+...$ and $h_{q+1}=id_{q+1}+t{\xi}_{q+1}+ ...$ when $t\rightarrow 0$, we obtain the local formula:ÊÊ\\
\[      d_{\mu}{\eta}^k={\eta}^k_rf^r_{\mu}+ ... ={\xi}^i({\partial}_if^k_{\mu}-f^k_{\mu+1_i})+f^k_{\mu+1_r}{\xi}^r+ ...+f^k_r{\xi}^r_{\mu}    \]
and thus ${\eta}_{q+1}=f_{q+2}({\bar{\xi}}_{q+1})$ with ${\bar{\xi}}_{q+1}={\xi}_{q+1}+{\chi}_{q+1}(\xi)$. This transformation is invertible if and only if $\xi \rightarrow {\bar{\xi}}=\xi+{\chi}_0(\xi)=A(\xi)$ is an 
isomorphism  of $T$.  \hspace{1cm}  $\Box$ \\

\noindent
{\bf EXAMPLE} 6.15: For $q=1$, we otain:  \\
\[ \begin{array}{ll}
\delta{\chi}^k_{,i}& =({\partial}_i{\xi}^k-{\xi}^k_i)+({\xi}^r{\partial}_r{\chi}^k_{,i}+{\chi}^k_{,r}{\partial}_i{\xi}^r-{\chi}^r_{,i}{\xi}^k_r)\\
  &=({\partial}_i{\bar{\xi}}^k-{\bar{\xi}}^k_i)+({\chi}^k_{r,i}{\bar{\xi}}^r-{\chi}^r_{,i}{\bar{\xi}}^k_r)     \\
\delta {\chi}^k_{j,i}&=({\partial}_i{\xi}^k_j-{\xi}^k_{ij})+({\xi}^r{\partial}_r{\chi}^k_{j,i}+{\chi}^k_{j,r}{\partial}_i{\xi}^r+{\chi}^k_{r,i}{\xi}^r_j-{\chi}^r_{j,i}{\xi}^k_r-{\chi}^r_{,i}{\xi}^k_{jr})   \\
  & =({\partial}_i{\bar{\xi}}^k_j-{\bar{\xi}}^k_{ij})+({\chi}^k_{rj,i}{\bar{\xi}}^r+{\chi}^k_{r,i}{\bar{\xi}}^r_j-{\chi}^r_{j,i}{\bar{\xi}}^k_r-{\chi}^r_{,i}{\bar{\xi}}^k_{jr})   
  \end{array}  \]
For the Killing system $R_1\subset J_1(T)$ with $g_2=0$, these variations are {\it exactly} the ones that can be found in ([C], (50)+(49), p 124 with a printing mistake corrected on p 128) when replacing a $3\times 3$ skewsymmetric matrix by the corresponding vector. {\it The last unavoidable Proposition is thus essential in order to bring back the nonlinear framework of finite elasticity to the linear framewok of infinitesimal elasticity that only depends on the linear Spencer operator}.\\
For the conformal Killing system ${\hat{R}}_1\subset J_1(T)$ (see next section) we obtain:    \\
\[ {\alpha}_i={\chi}^r _{r,i}\Rightarrow \delta {\alpha}_i=({\partial}_i{\xi}^r_r-{\xi}^r_{ri})+({\xi}^r{\partial}_r{\alpha}_i+{\alpha}_r{\partial}_i{\xi}^r+{\chi}^s_{,i}{\xi}^r_{rs})  \]
This is {\it exactly} the variation obtained by Weyl ([W], (76), p 289) who was assuming implicitly $A=0$ when setting ${\bar{\xi}}^r_r=0\Leftrightarrow {\xi}^r_r=-{\alpha}_i{\xi}^i$ by introducing a connection. Accordingly, ${\xi}^r_{ri}$ is the variation of the EM potential itself, that is the $\delta A_i$ of engineers used in order to exhibit the Maxwell equations from a variational principle ([W], $\S$ 26) but the introduction of the Spencer operator is new in this framework.\\

Finally, chasing in diagram $(1)$ , the Spencer sequence is locally exact at $C_1$ if and only if the Janet sequence is locally exact at $F_0$ because the central sequence is locally exact. {\it The situation is much more complicate in the nonlinear framewok}. Let $\bar{\omega}$ be a section of $\cal{F}$ satisfying the same CC as $\omega$, namely $I(j_1(\omega))=0$. It follows that we can find a section $f_{q+1}\in {\Pi}_{q+1}$ such that $f^{-1}_q(\omega)=\bar{\omega}\Rightarrow j_1(f^{-1}_q)(j_1(\omega))=j_1(f^{-1}_q(\omega))=j_1({\bar{\omega}})$ and $f^{-1}_{q+1}(j_1(\omega))=j_1(\bar{\omega})$. We obtain therefore 
$ j_1(f^{-1}_q)(j_1(\omega))=f^{-1}_{q+1}(j_1(\omega)) \Rightarrow (f_{q+1}\circ  j_1(f^{-1}_q))^{-1}(j_1(\omega))-j_1(\omega)=L({\sigma}_q)\omega=0$ and thus ${\sigma}_q={\bar{D}}f^{-1}_{q+1}\in T^*\otimes R_q $ {\it over the target}, even if $f_{q+1} $ may not be a section of ${\cal{R}}_{q+1}$. As ${\sigma}_q$ is killed by ${\bar{D}}'$, we have related cocycles at $\cal{F}$ in the Janet sequence {\it over the source} with cocycles at $T^*\otimes R_q$ or 
$C_1$ {\it over the target}.\\
Now, if $f_{q+1},f'_{q+1}\in {\Pi}_{q+1}$ are such that $f^{-1}_{q+1}(j_1(\omega))=f'^{-1}_{q+1}(j_1(\omega))=j_1({\bar{\omega}})$, it follows that $(f'_{q+1}\circ f^{-1}_{q+1})(j_1(\omega))=j_1(\omega)   \Rightarrow
\exists  g_{q+1}\in {\cal{R}}_{q+1}$ such that $f'_{q+1}=g_{q+1}\circ f_{q+1}$ and the new ${\sigma}'_q=\bar{D}f'^{-1}_{q+1}$ differs from the initial ${\sigma}_q=\bar{D}f^{-1}_{q+1}$ by a gauge transformation.\\
Conversely, let $f_{q+1},f'_{q+1}\in {\Pi}_{q+1}$ be such that ${\sigma}_q=\bar{D}f^{-1}_{q+1}=
\bar{D}f'^{-1}_{q+1}={\sigma}'_q$. It follows that $\bar{D}(f^{-1}_{q+1}\circ f'_{q+1})=0$ and one can find $g\in aut(X)$ such that $f'_{q+1}=f_{q+1}\circ j_{q+1}(g)$ providing ${\bar{\omega}}'=f'^{-1}_q(\omega)=(f_q\circ j_q(g))^{-1}(\omega)=j_q(g)^{-1}(f^{-1}_q(\omega))=j_q(g)^{-1}(\bar{\omega})$.\\

\noindent
{\bf PROPOSITION} 6.16: Natural transformations of $\cal{F}$ {\it over the source} in the nonlinear Janet sequence correspond to gauge transformations of $T^*\otimes R_q$ or $C_1$ {\it over the target} in the nonlinear Spencer sequence. Similarly, the Lie derivative ${\cal{D}}\xi={\cal{L}}(\xi)\omega\in F_0$ in the linear Janet sequence corresponds to the Spencer operator $D{\xi}_{q+1}\in T^*\otimes R_q$ or $D_1{\xi}_q\in C_1$ in the linear Spencer sequence.\\

With a slight abuse of language $\delta f=\eta\circ f\Leftrightarrow \delta f\circ f^{-1}=\eta \Leftrightarrow f^{-1}\circ \delta f=\xi$ when $\eta=T(f)(\xi)$ and we get $j_q(f)^{-1}(\omega)=\bar{\omega} \Rightarrow j_q(f+\delta f)^{-1}(\omega)=\bar{\omega}+\delta \bar{\omega}$ that is $j_q(f^{-1}\circ(f+\delta f))^{-1}(\bar{\omega})=\bar{\omega}+\delta \bar{\omega} \Rightarrow \delta \bar{\omega}={\cal{L}}(\xi)\bar{\omega}$ and $j_q((f+\delta f)\circ f^{-1}\circ f)^{-1}(\omega)=j_q(f)^{-1}(j_q((f+\delta f)\circ f^{-1})^{-1}(\omega)) \Rightarrow \delta \bar{\omega}=j_q(f)^{-1}({\cal{L}}(\eta)\omega)$.\\
Passing to the infinitesimal point of view, we obtain the following generalization of Remark 3.12 which is important for applications ([2], AJSE-mathematics):\\

\noindent
{\bf COROLLARY} 6.17: $\delta \bar{\omega}={\cal{L}}(\xi)\bar{\omega}=j_q(f)^{-1}({\cal{L}}(\eta)
\omega)$.\\

\noindent
{\bf EXAMPLE} 6.18: In Example 3.4 with $n=1, q=1$, we have $\omega(f(x))f_x(x)={\bar{\omega}}(x), \omega(f(x))f_{xx}(x)+ {\partial}_y\omega(f(x))f^2_x(x)={\partial}_x{\bar{\omega}}(x)$ and obtain 
therefore $ \omega {\sigma}_{y,y}+{\sigma}_{,y}{\partial}_y\omega \equiv -\omega(1/f_x)({\partial}_xf_x-f_{xx})(1/{\partial}_xf)+((f_x/{\partial}_xf)-1){\partial}_y\omega=0$ whenever $y=f(x)$. The case of an affine stucture needs more work.Ê \\

\noindent
{\bf 7   COSSERAT VERSUS WEYL: NEW PERSPECTIVES FOR PHYSICS}\\

   As an application of the previous mehods, let us now consider the {\it conformal Killing system}:\\ 
 \[  \begin{array}{ll}
   {\hat{R}}_1\subset J_1(T)  \hspace{4cm} & {\omega}_{rj}{\xi}^r_i+{\omega}_{ir}{\xi}^r_j+{\xi}^r{\partial}_r{\omega}_{ij}=A(x){\omega}_{ij}  \hspace{2cm}
   \end{array}  \]
 with symbols:\\
 \[  {\hat{g}}_2\subset S_2T^*\otimes T  \hspace {2cm} n{\xi}^k_{ij}={\delta}^k_i{\xi}^r_{rj}+{\delta}^k_j{\xi}^r_{ri}-{\omega}_{ij}{\omega}^{ks}{\xi}^r_{rs}\hspace{5mm} \Rightarrow \hspace{5mm}{\hat{g}}_3=0, \forall n\geq 3 \]
 obtained by eliminating the arbitrary function $A(x)$, where $\omega$ is the Euclidean metric when $n=1$ (line), $n=2$ (plane) or $n=3$ (space) and the Minskowskian metric when $n=4$ (space-time).\\
  The brothers Cosserat were only dealing with the {\it Killing subsystem}:  \\
  \[ R_1\subset {\hat{R}}_1 \hspace{5cm} {\omega}_{rj}{\xi}^r_i+{\omega}_{ir}{\xi}^r_j+{\xi}^r{\partial}_r{\omega}_{ij}=0\hspace{3cm} \]
 that is with $\{{\xi}^k,{\xi}^k_i\mid {\xi}^r_r=0, {\xi}^k_{ij}=0\}=\{translations, rotations\}$ when $A(x)=0$, while, {\it in a somehow complementary way}, Weyl was mainly dealing with $\{{\xi}^r_r, {\xi}^r_{ri}\}=\{dilatation, elations\}$. Accordingly, one has ([7]):\\
 
 \noindent
 {\bf THEOREM} 7.1: The {\it Cosserat equations} ([C], p 137 for $n=3$, p 167 for $n=4$): \\
 \[ {\partial}_r{\sigma}^{i,r}=f^i  \hspace{1cm}, \hspace{1cm}{\partial}_r{\mu}^{ij,r}+{\sigma}^{i,j}-{\sigma}^{j,i}=m^{ij} \]
 are {\it exactly} described by the formal adjoint of the first Spencer operator $D_1:R_1\rightarrow T^*\otimes R_1$. Introducing ${\phi}^{r,ij}=-{\phi}^{r,ji}$ and ${\psi}^{rs,ij}=-{\psi}^{rs,ji}=-{\psi}^{sr,ij}$, they can be {\it parametrized} by the formal adjoint of the second Spencer operator $D_2:T^*\otimes R_1\rightarrow {\wedge}^2T^*\otimes R_1$:\\
 \[   {\sigma}^{i,j}={\partial}_r{\phi}^{i,jr}  \hspace{1cm}, \hspace{1cm} {\mu}^{ij,r}=
 {\partial}_s{\psi}^{ij,rs}+{\phi}^{j,ir}-{\phi}^{i,jr}   \]
 
 \noindent
 {\bf EXAMPLE} 7.2: When $n=2$, lowering the indices by means of the constant metric $\omega$, we just need to look for the factors of ${\xi}_1,{\xi}_2$ and ${\xi}_{1,2}$ in the integration by parts of the sum:\\
 \[ {\sigma}^{1,1}({\partial}_1{\xi}_1-{\xi}_{1,1})+{\sigma}^{1,2}({\partial}_2{\xi}_1-{\xi}_{1,2})+{\sigma}^{2,1}({\partial}_1{\xi}_2-{\xi}_{2,1})+{\sigma}^{2,2}({\partial}_2{\xi}_2-{\xi}_{2,2})+{\mu}^{12,r}({\partial}_r{\xi}_{1,2}-{\xi}_{1,2r}) \]
 Finally, setting ${\phi}^{1,12}={\phi}^1,{\phi}^{2,12}={\phi}^2,{\psi}^{12,12}={\phi}^3$, we obtain the nontrivial parametrization ${\sigma}^{1,1}={\partial}_2{\phi}^1, {\sigma}^{1,2}=-{\partial}_1{\phi}^1, {\sigma}^{2,1}=-{\partial}_2{\phi}^2, {\sigma}^{2,2}={\partial}_1{\phi}^2, {\mu}^{12,1}={\partial}_2{\phi}^3+{\phi}^1, {\mu}^{12,2}=-{\partial}_1{\phi}^3-{\phi}^2$ in a coherent way with the Airy parametrization obtained when ${\phi}^1={\partial}_2{\phi}, {\phi}^2={\partial}_1{\phi}, {\phi}^3=-\phi$.\\

\noindent
{\bf REMARK} 7.3 : First of all, it is clear that [C] (p 13,14 for $n=1$, p 75,76 for $n=2$) still deals with $m=3$ for the "ambient space", that is with the construction of the nonlinear gauge sequence, in particular for the dynamical study of a line with arc length $s$ and time $t$ considered as a surface, hence with no way to pass from the source to the target, only possible, as we have seen, when $m=n=3$ by using the nonlinear Spencer sequence. For $n=4$, the group of rigid motions of space is extended to space-time by using only a translation on time and we can rewrite the formulas in ([C], p 167) as follows:  \\
\[  \frac{d}{dt}=\frac{dx}{dt}\frac{\partial}{\partial x}+\frac{dy}{dt}\frac{\partial}{\partial y}+\frac{dz}{dt}\frac{\partial}{\partial z}+\frac{\partial}{\partial t} \Rightarrow
 \frac{\partial p_{xx}}{\partial x}+...+\frac{1}{\Delta}\frac{dA}{dt}=\frac{\partial}{\partial x}(p_{xx}+\frac{A}{\Delta}\frac{dx}{dt})+...+\frac{\partial}{\partial t}(\frac{A}{\Delta})\]
 \[  \frac{\partial q_{xx}}{\partial x}+...+p_{yz}-p_{zy}+\frac{1}{\Delta}\frac{dP}{dt}+\frac{C}{\Delta}\frac{dy}{dt}-\frac{B}{\Delta}\frac{dz}{dt}=\frac{\partial}{\partial x}(q_{xx}+\frac{P}{\Delta}\frac{dx}{dt})+...+\frac{\partial}{\partial t}(\frac{P}{\Delta})+(p_{yz}+\frac{C}{\Delta}\frac{dy}{dt})-(p_{zy}+\frac{B}{\Delta}\frac{dz}{dt})  \]
 It is essential to notice that the Cosserat equations for $n=3$ are still introduced today in a phenomenological way ([35] is a good example), contrary to the "deductive" way used in ([C], p 1-6) and that "intuition" will never allow to provide the relativistic Cosserat equations for $n=4$ which are presented for the first time.\\

\noindent
{\bf THEOREM} 7.4: The {\it Weyl equations} ([W], \S 35) are {\it exactly} described by the formal adjoint of the first Spencer operator $D_1:{\hat{R}}_2\rightarrow T^*\otimes {\hat{R}}_2$ when $n=4$ and can be parametrized by the formal adjoint of the second Spencer operator $D_2:T^*\otimes {\hat{R}}_2\rightarrow T^*\otimes {\hat{R}}_2$. In particular, among the components of the first Spencer operator, one has ${\partial}_i{\xi}^r_{rj}-{\xi}^r_{ijr}={\partial}_i{\xi}^r_{rj}$ and thus the components ${\partial}_i{\xi}^r_{rj}-{\partial}_j{\xi}^r_{ri}=F_{ij}$ of the EM {\it field} with EM {\it potential} ${\xi}^r_{ri}=A_i$ coming from the second order jets ({\it elations}). It follows that $D_1$ projects onto $d:T^*\rightarrow {\wedge}^2T^*$ and thus $D_2$ {\it projects onto} the first set of Maxwell equations described by $d:{\wedge}^2T^*\rightarrow {\wedge}^3T^*$. Indeed, the Spencer sequence projects onto the Poincar\'{e} sequence {\it with a shift by $+1$ in the degree of the exterior forms involved} because both sequences are made with first order involutive operators and the comment after diagram (1) can thus be used. By duality, the second set of Maxwell equations thus {\it appears among} the Weyl equations which {\it project onto} the Cosserat equations because of the inclusion $R_1\simeq R_2\subset {\hat{R}}_2$.\\

\noindent
{\bf REMARK} 7.5: When $n=4$, the {\it Poincar\'{e} group} ($10$ parameters) is a subgroup of the {\it conformal group} ($15$ parameters) which is not a maximal subgroup because it is a subgroup of the {\it Weyl group} ($11$ parameters) obtained by adding the only dilatation with infinitesimal generators $x^i{\partial}_i$. However, the {\it optical group} is another subgroup with $10$ parameters which is maximal and the same procedure may be applied to all these subgroups in order to study coupling phenomena. It is also important to notice that {\it the first and second sets of Maxwell equations are invariant by any diffeomorphism and the conformal group is only the group of invariance of the Minkowski constitutive laws in vacuum} ([20])([27], p 492).\\

\noindent
{\bf REMARK} 7.6: Though striking it may look like, {\it there is no conceptual difference between the Cosserat and Maxwell equations on space-time}. As a byproduct, separating space from time, there is no conceptual difference between the Lam\'e constants (mass per unit volume) of elasticity and the magnetic (dielectric) constants of EM appearing in the respective wave speeds. For example, the speed of longitudinal free vibrations of a thin elastic bar with Young modulus $E$ and mass per unit volume $\rho$ is $v=\sqrt{\frac{E}{\rho}}$ while the speed of light in a medium with magnetic constant $\mu$ and dielectric constant $\epsilon$ is $v=\sqrt{\frac{1/\mu}{\epsilon}}$. In the first case, we have the $1$-dimensional dynamical equations: \\
\[      \delta \int (\frac{1}{2}E(\frac{\partial \xi}{\partial x})^2-\frac{1}{2}\rho(\frac{\partial \xi}{\partial t})^2)dxdt=0  \Rightarrow  E\frac{{\partial}^2\xi}{\partial x^2}-\rho \frac{{\partial}^2\xi}{\partial t^2}=0  \]
In the second case, studying the propagation in vacuum for simplicity, one uses to set $\vec{H}=(1/{\mu}_0)\vec{B}, \vec{D}={\epsilon}_0 \vec{E}$ with ${\epsilon}_0 {\mu}_0c^2=1$ in the induction equations and to substitute the space-time parametrization $dA=F$ of the field equations $dF=0$ in the variational condition $\delta \int (\frac{1}{2}{{\epsilon}_0}{\vec{E}}^2-\frac{1}{2}(1/{\mu}_0){\vec{B}}^2)dxdt=0$. However, the second order PD equations thus obtained become wave equations {\it only if one assumes the Lorentz condition} $div(A)={\omega}^{ij}{\partial}_iA_j=0$ ([20]). This is not correct because the Lagrangian of the corresponding {\it variational problem with constraint} must contain the additional term $\lambda div(A)$ where $\lambda $ is a {\it Lagrange multiplier} providing the equations $\Box A=d\lambda $ as a $1$-form and thus $ \Box F=0$ as a $2$-form when $\Box$ is the Dalembertian ([27], p 885).\\

\noindent
{\bf REMARK} 7.7: When studying static phenomena, $\epsilon=({\epsilon}_{ij})$ and $\vec{E}=(E^i)$ are now on equal footing in the Lagrangian, exactly like in the technique of finite elements. Starting with a homogeneous medium at rest with no stress and electric induction, we may consider a quadratic Lagrangian $A^{ijkl}{\epsilon}_{ij}{\epsilon}_{kl}+B^{ij}E_iE_j+C^{ijk}{\epsilon}_{ij}E_k$ obtained by moving the indices by means of the Euclidean metric. The two first terms describe (pure) linear elasticity and electrostatic while only the last {\it quadratic coupling term} may be used in order to describe coupling phenomena. For an isotropic medium, the $3$-tensor $C$ must vanish and such a coupling phenomenon, called {\it piezzoelectricity}, can only appear in non-isotropic media like crystals, providing the additional stress ${\sigma}^{ij}=C^{ijk}E_k$ and/or an additional electric induction $D^k=C^{ijk}{\epsilon}_{ij}$. Accordingly, if the medium is fixed, for example between the plates of a condenser, an electric field may provide stress inside while, if the medium is deformed as in the piezzo-lighters, an electric induction may appear and produce a spark. Finally, for an isotropic medium, we can only add a {\it cubic coupling term} $C^{ijkl}{\epsilon}_{ij}E_kE_l$ responsible for {\it photoelasticity} as it provides the additional electric induction $D^l=(C^{ijkl}{\epsilon}_{ij})E_k$, modifying therefore the dielectric constant by a term depending linearly on the deformation and thus modifying the index of refraction $n$ because $\epsilon{\mu}_0c^2=n^2$ with ${\epsilon}_0{\mu}_0c^2=1$ in vacuum leads to $\epsilon=n^2{\epsilon}_0$. We may also identify the dimensionless "speed" $v^k/c\ll 1, \forall k=1,2,3$ (time derivative of position) with a first jet (Lorentz rotation) by setting ${\partial}_4{\xi}^k-{\xi}^k_4=0$ and introduce the {\it speed of deformation} by the formula $2{\nu}_{ij}={\omega}_{rj}({\partial}_i{\xi}^r_4-{\xi}^r_{i4})+{\omega}_{ir}({\partial}_j{\xi}^r_4-{\xi}^r_{j4})={\omega}_{rj}{\partial}_i{\xi}^r_4+{\omega}_{ir}{\partial}_j{\xi}^r_4={\partial}_4({\omega}_{rj}{\partial}_i{\xi}^r+{\omega}_{ir}{\partial}_j{\xi}^r)={\omega}_{rj}{\partial}_iv^r+{\omega}_{ir}{\partial}_jv^r={\partial}_4{\epsilon}_{ij}, \forall 1\leq i,j\leq 3$ in order to obtain {\it streaming birefringence} in a similar way. These results perfectly agree with most of the field-matter couplings known in engineering sciences ([28]) but contradict {\it gauge theory} ([15],[26]) and {\it general relativity} ([W],[21]).\\

In order to justify the last remark, let $G$ be a Lie group with identity $e$ and parameters $a$ acting on $X$ through the group action $X\times G\rightarrow X: (x,a)\rightarrow y=f(x,a)$ and (local) infinitesimal generators ${\theta}_{\tau}$ satisfying $[{\theta}_{\rho},{\theta}_{\sigma}]=c^{\tau}_{\rho \sigma}{\theta}_{\tau}$ for $\rho,\sigma,\tau =1,...,dim(G)$. We may prolong the {\it graph} of this action by differentiating $q$ times the action law in order to eliminate the parameters in the following commutative and exact diagram where ${\cal{R}}_q$ is a Lie groupoid with local coordinates $(x,y_q)$, {\it source} projection ${\alpha}_q:(x,y_q)\rightarrow (x)$ and {\it target} projection ${\beta}_q:(x,y_q)\rightarrow (y)$ when $q$ is large enough:\\

\[  \begin{array}{rclcccl}
0\rightarrow & X\times G  &  \longrightarrow         &   &  {\cal{R}}_q &  & \rightarrow 0 \\
                      &  \parallel    &    &  {\alpha}_q\swarrow  &   & \searrow {\beta}_q  &     \\
    &   X   \times G  & \rightarrow  &  X  & \times & X    
    \end{array}  \]

 The link between the various sections of the trivial principal bundle on the left ({\it gauging procedure}) and the various corresponding sections of the Lie groupoid on the right with respect to the source projection is expressed by the next commutative and exact diagram:\\
  
\[  \begin{array}{rclccl}
0\rightarrow  & \hspace{13mm} X  \hspace{6mm} \times & G &  =  &  \hspace{4mm}{\cal{R}}_q  &  \rightarrow 0  \\
      &   a=cst \uparrow \downarrow \uparrow a(x)  &  &  &  j_q(f) \uparrow \downarrow \uparrow  f_q &  \\
         &  \hspace{4mm}X &    &  =   &    \hspace{4mm}X   &   
\end{array}  \]

\noindent
{\bf THEOREM} 7.8: In the above situation, the nonlinear Spencer sequence is isomorphic to the nonlinear gauge sequence and we have the following commutative and locally exact diagram:\\
\[ \begin{array}{rccccc}
  &   X\times G & \rightarrow & T^*\otimes {\cal{G}} &\stackrel{MC}{ \rightarrow} &  {\wedge}^2T^*\otimes {\cal{G}}  \\
    & \downarrow  &   &  \downarrow  &  & \downarrow  \\
    0\rightarrow \Gamma \rightarrow & {\cal{R}}_q  & \stackrel{\bar{D}}{\rightarrow} & T^*\otimes R_q & \stackrel{{\bar{D}}'}{\rightarrow}  & {\wedge}^2T^* \otimes R_q  
    \end{array}     \]
{\it The action is essential in the Spencer sequence but disappears in the gauge sequence}.\\

\noindent
{\bf Proof}: If we consider the action $y=f(x,a)$ and start with a section $(x)\rightarrow (x,a(x))$ of $X\times G$, we obtain the section $(x)\rightarrow (x,f^k_{\mu}(x)={\partial}_{\mu}f^k(x,a(x)))$ of ${\cal{R}}_q$. Setting $b=a^{-1}=b(a)$, we get $y=f(x,a)\Rightarrow x=f(y,b)\Rightarrow y=f(f(y,b(a),a)$ and thus $\frac{\partial y}{\partial x}\frac{\partial f}{\partial b}\frac{\partial b}{\partial a}+\frac{\partial y}{\partial a}=0$ with $\frac{\partial f}{\partial b}=\theta(x)\omega(b) $ from the first fundamental theorem of Lie. With $-\omega(b)db=-dbb^{-1}=a^{-1}da$, we obtain:   \\
\[ \begin{array}{rl}
  {\partial}_if^k_{\mu}-f^k_{\mu +1_i}  &  = d_i({\partial}_{\mu}f^k(x,a(x))-{\partial}_{\mu+1_i}f^k(x,a(x)) \\
           & = {\partial}_{\mu}(\frac{\partial f^k}{\partial a^{\tau}}){\partial}_ia^{\tau}  \\
           &=-{\partial}_{\mu}(\frac{\partial f^k}{\partial x^r}{\theta}^r_{\tau}(x)){\omega}^{\tau}_{\sigma}(b)\frac{\partial b^{\sigma}}{\partial a^{\tau}}{\partial}_ia^{\tau}
    \end{array}     \]
and thus ${\chi}^k_{\mu,i}(x)=A^{\tau}_i(x){\partial}_{\mu}{\theta}^k_{\tau}(x)$ from the inductive formula allowing to define ${\chi}_q=\bar{D}f_{q+1}$. \\
As for the commutatitvity of the right square, we have: \\
\[ {\partial}_i{\chi}^k_{\mu ,j}-{\partial}_j{\chi}^k_{\mu ,i}-{\chi}^k_{\mu +1_i,j}+{\chi}^k_{\mu +1_j,i}=({\partial}_i{A}^{\tau}_j-{\partial}_j{A}^{\tau}_i){\partial}_{\mu}{\theta}^k_{\tau} \]
\[  (\{{\chi}_{q+1}({\partial}_i),{\chi}_{q+1}({\partial}_j)\})^k_{\mu}={A}^{\rho}_i{A}^{\sigma}_j{\partial}_{\mu}([{\theta}_{\rho},{\theta}_{\sigma}])^k=c^{\tau}_{\rho\sigma}{A}^{\rho}_i{A}^{\sigma}_j{\partial}_{\mu}{\theta}^k_{\tau} . \hspace{1cm} \Box   \]

Introducing now the Lie algebra ${\cal{G}}=T_e(G)$ and the Lie algebroid $R_q\subset J_q(T)$, namely the linearization of ${\cal{R}}_q$ at the $q$-jet of the identity $y=x$, we get the commutative and exact diagram:\\

\[  \begin{array}{rclccl}
0\rightarrow  &  \hspace{13mm} X   \hspace{6mm} \times & \cal{G} & = & \hspace{4mm} R_q & 
\rightarrow 0 \\
  &  \lambda=cst \uparrow \downarrow \uparrow \lambda(x) & & & j_q(\xi) \uparrow \downarrow \uparrow {\xi}_q &  \\
    &  \hspace{4mm} X &  & = & \hspace{4mm} X  & 
    \end{array}  \]
\noindent
where the upper isomorphism is described by ${\lambda}^{\tau}(x)\rightarrow {\xi}^k_{\mu}(x)={\lambda}^{\tau}(x){\partial}_{\mu}{\theta}^k_{\tau}(x)$ for $q$ large enough. The unusual Lie algebroid structure on $X\times {\cal{G}}$ is described by the formula: $([{\lambda},{\lambda}'])^{\tau}=c^{\tau}_{\rho \sigma}{\lambda}^{\rho}{\lambda}'^{\sigma} + ({\lambda}^{\rho}{\theta}_{\rho}).{\lambda}'^{\tau}-({\lambda}'^{\sigma}{\theta}_{\sigma}).{\lambda}^{\tau}$ which is induced by the ordinary bracket $[{\xi},{\xi}']$ on $T$ and thus depends on the action. Applying the Spencer operator, we finally obtain ${\partial}_i{\xi}^k_{\mu}(x)-{\xi}^k_{\mu +1_i}(x)={\partial}_i{\lambda}^{\tau}(x){\partial}_{\mu}{\theta}^k_{\tau}(x)$ and the linear Spencer sequence is isomorphic to the linear gauge sequence already introduced which is no longer depending on the action as it is only the tensor product of the Poincar\'{e} sequence by $\cal{G}$.  \\

\noindent
{\bf EXAMPLE} 7.9: Let us consider the group of affine transformations of the real line $y=a^1x+a^2$ with $n=1, dim(G)=2, q=2$, ${\cal{R}}_2$ defined by the system $y_{xx}=0$, $R_2$ defined by ${\xi}_{xx}=0$ and the two infinitesimal generators ${\theta}_1=x\frac{\partial}{\partial x}, {\theta}_2=\frac{\partial}{\partial x}$. We get $f(x)=a^1(x)x+a^2(x), f_x(x)=a^1(x), f_{xx}(x)=0$ and thus ${\chi}_{,x}(x)=(1/f_x(x)){\partial}_xf(x)-1=(1/a^1(x))(x{\partial}_xa^1(x)+{\partial}_xa^2(x))=xA^1_x(x)+A^2_x(x), {\chi}_{x,x}(x)=(1/f_x(x))({\partial}_xf_x(x)-(1/f_x(x)){\partial}_xf(x)f_{xx}(x))=(1/a^1(x)){\partial}_xa^1(x)=A^1_x(x), 
{\chi}_{xx,x}(x)=0$. Similarly, we get $\xi(x)={\lambda}^1(x)x+{\lambda}^2(x), {\xi}_x(x)={\lambda}^1(x), {\xi}_{xx}(x)=0$. Finally, integrating by part the sum $\sigma ({\partial}_x\xi-{\xi}_x)+\mu({\partial}_x{\xi}_x-{\xi}_{xx})$ we obtain the dual of the Spencer operator as ${\partial}_x\sigma=f,{\partial}_x\mu+\sigma=m$ that is to say the Cosserat equations for the affine group of the real line.\\

     It finally remains to study GR within this framework, as it is only "{\it added}" by Weyl in an independent way and, for simplicity, we shall restrict to the linearized aspect. First of all, it becomes clear from diagram (1) that the mathematical foundation of GR is based on a confusion between the operator ${\cal{D}}_1$ ({\it classical curvature alone}) in the Janet sequence when $\cal{D}$ is the Killing operator brought to involution and the operator $D_2$ ({\it gauge curvature=curvature+torsion}) in the corresponding Spencer sequence. It must also be noticed that, according to the same diagram, the bigger is the underlying group, the bigger are the Spencer bundles while, on the contrary, the smaller are the Janet bundles depending on the invariants of the group action (deformation tensor in classical elasticity is a good example). Precisely, as already noticed in Theorem 7, if $G\subset \hat{G}$, the Spencer sequence for $G$ is {\it contained into} the Spencer sequence for $\hat{G}$ while the Janet sequence for $G$ {\it projects onto} the Janet sequence for $\hat{G}$, {\it the best picture for understanding such a phenomenon is that of two children sitting on the ends of a beam and playing at see-saw}.\\
      
      Such a confusion is also combined with another one well described in ([40], p 631) by the chinese saying "{\it To put Chang's cap on Li's head}", namely to relate the Ricci tensor (usually obtained from the Riemann tensor by contraction of indices) to the energy-momentum tensor (space-time stress), without taking into account the previous confusion relating the gauge curvature to {\it rotations} only while the (classical and Cosserat) stress has only to do with {\it translations}. In addition, it must be noticed that {\it the Cosserat and Maxwell equations can be parametrized while the Einstein equations cannot be parametrized} ([29]).\\

        In order to escape from this dilemna, let us denote by $B^2(g_q), Z^2(g_q)$ and $H^2(g_q)=Z^2(g_q)/B^2(g_q)$ the coboundary (image of the left $\delta$), cocycle (kernel of the right $\delta$) and cohomology bundles of the $\delta$-sequence  $T^*\otimes g_{q+1} \stackrel{\delta}{\rightarrow} {\wedge}^2T^*\otimes g_q \stackrel{\delta}{\rightarrow} {\wedge}^3T^*\otimes S_{q-1}T^*\otimes T$. It can be proved that the order of the generating CC of a formally integrable operator of order $q$ is equal to $s+1$ when $s$ is the smallest integer such that $H^2(g_{q+r})=0, \forall r\geq s$ ([26]). As an example with $n=3$, we let the reader prove that the second order systems $y_{33}=0, y_{23}-y_{11}=0, y_{22}=0$ and $y_{33}-y_{11}=0, y_{23}=0, y_{22}-y_{11}=0$ have both three second order generating CC ([30]). For the Killing  system $R_1\subset J_1(T)$ with symbol $g_1$, we have $F_0=J_1(T)/R_1=T^*\otimes T/g_1$ and the short exact sequence $0\rightarrow g_1\rightarrow T^*\otimes T\rightarrow F_0\rightarrow 0$. As $q=1$ and $g_2=0\Rightarrow g_3=0$ we have $s=1$ and no CC of order $1$. The generating CC of order $2$ only depend on $F_1={\omega}^{-1}({\cal{F}}_1)$ according to section 2 where $F_1$ is now defined by the following commutative diagram with exact columns but the first on the left and exact rows:\\
 
\[  \begin{array}{rcccccccl}
   &  0 & & 0 & & 0 &  &  &   \\
   & \downarrow & & \downarrow & & \downarrow & & &  \\
0\rightarrow & g_3 & \rightarrow &  S_3T^*\otimes T & \rightarrow & S_2T^*\otimes F_0& \rightarrow & F_1 & \rightarrow 0  \\
   & \hspace{2mm}\downarrow  \delta  & & \hspace{2mm}\downarrow \delta & &\hspace{2mm} \downarrow \delta & & &  \\
0\rightarrow& T^*\otimes g_2&\rightarrow &T^*\otimes S_2T^*\otimes T & \rightarrow &T^*\otimes T^*\otimes F_0 &\rightarrow & 0 &  \\
   &\hspace{2mm} \downarrow \delta &  &\hspace{2mm} \downarrow \delta & &\hspace{2mm}\downarrow \delta &  &  &   \\
0\rightarrow & {\wedge}^2T^*\otimes g_1 & \rightarrow & {\wedge}^2T^*\otimes T^*\otimes T & \rightarrow & {\wedge}^2T^*\otimes F_0 & \rightarrow & 0 &  \\
   &\hspace{2mm}\downarrow \delta  &  & \hspace{2mm} \downarrow \delta  &  & \downarrow  & &  &  \\
0\rightarrow & {\wedge}^3T^*\otimes T & =  & {\wedge}^3T^*\otimes T  &\rightarrow   & 0  &  &  &   \\
    &  \downarrow  &  &  \downarrow  &  &  &  &  &  \\
    &  0  &   & 0  & &  &  &  &
\end{array}  \]

It follows from a chase([26], p 55)([27], p 192)([32], p 171) that there is a short exact {\it connecting sequence} $0\rightarrow B^2(g_1) \rightarrow Z^2(g_1) \rightarrow F_1  \rightarrow 0$ leading to an isomorphism $F_1\simeq H^2(g_1)$. The Riemann tensor is thus a section of $Riemann=F_1=H^2(g_1)=Z^2(g_1)$ in the Killing case with $dim(Riemann)=(n^2(n+1)^2/4)-(n^2(n+1)(n+2)/6)=(n^2(n-1)^2/4)-(n^2(n-1)(n-2)/6)=n^2(n^2-1)/12$ by using either the upper row or the left column and {\it we find back the two algebraic properties of the Riemann tensor without using indices}. \\
However, for the conformal Killing system, we still have $q=1$ but the situation is much more delicate because $g_3=0$ for $n\geq 3$ and $H^2({\hat{g}}_2)=0$ only for $n\geq 4$ ([26], p 435). Hence, setting similarly ${\hat{F}}_0=T^*\otimes T/{\hat{g}}_1$, the Weyl tensor is a section of $Weyl={\hat{F}}_1=H^2({\hat{g}}_1)\neq Z^2({\hat{g}}_1)$. The inclusion $g_1\subset {\hat{g}}_1$ and the relations $g_2=0, {\hat{g}}_3=0$ finally induce the following {\it crucial} commutative and exact {\it diagram} (2) ([25], p 430):\\
              
\[ \begin{array}{rcccccccll}
 & & & & & & & 0 & &\\
 & & & & & & & \downarrow && \\
  & & & & & 0& & Ricci & & \\
  & & & & & \downarrow & & \downarrow  & & \\
   & & & 0 &\rightarrow & Z^2(g_1) & \rightarrow & Riemann  & \rightarrow 0 & \\
   & & & \downarrow & & \downarrow & &  \downarrow  &  &JANET\\
   & 0 &\rightarrow & T^*\otimes {\hat{g}}_2 & \stackrel{\delta}{\rightarrow} & Z^2({\hat{g}}_1) & \rightarrow & Weyl & \rightarrow 0 & \\
    & & & \downarrow & & \downarrow & & \downarrow &  & \\
 0 \rightarrow & S_2T^* & \stackrel{\delta}{\rightarrow}& T^*\otimes T^* &\stackrel{\delta}{\rightarrow} & {\wedge}^2T^* & \rightarrow & 0 &  & \\
   & & & \downarrow &  & \downarrow & & & & \\
   & & & 0 & & 0 & & &  &\\
   &&&& SPENCER &&&&&
   \end{array}  \]
   \noindent

A diagonal chase allows to identify $Ricci$ with $S_2T^*$ {\it without contracting indices} and provides the splitting of $T^*\otimes T^*$ into $S_2T^*$ ({\it gravitation}) and ${\wedge}^2T^*$ ({\it electromagnetism}) in the lower horizontal sequence {\it obtained by using the Spencer sequence}, solving thus an old conjecture. However, $T^*\otimes T^*\simeq T^*\otimes {\hat{g}}_2$ has only to do with second order jets ({\it elations}) and not a word is left from the standard approach to GR. In addition, we obtain the following important theorem explaining for the first time classical results in an intrinsic way:\\

\noindent
{\bf THEOREM} 7.10: There exist canonical splittings of the various $\delta$-maps appearing in the above diagram which allow to split the vertical short exact sequence on the right.\\

\noindent
{\bf Proof}: We recall first that a short exact sequence $0 \rightarrow M' \stackrel{f}{\rightarrow} M \stackrel{g}{\rightarrow} M" \rightarrow 0$ of modules {\it splits}, that is $M\simeq M'\oplus M"$, if and only if there exists a map $u:M \rightarrow M'$ with $u\circ f=id_{M'}$ or a map $v:M"\rightarrow M$ with $g\circ v=id_{M"}$ ([3] p 73)([32], p 33). Hence, starting with $({\tau}^k_{li,j})\in T^*\otimes {\hat{g}}_2$, we may introduce $({\rho}^k_{l,ij}={\tau}^k_{li,j}-{\tau}^k_{lj,i})\in B^2({\hat{g}}_1)\subset Z^2({\hat{g}}_1)\subset {\wedge}^2T^*\otimes {\hat{g}}_1$ but now ${\varphi}_{ij}={\rho}^r_{r,ij}={\tau}^r_{ri,j}-{\tau}^r_{rj,i}={\rho}_{ij}-{\rho}_{ji}\neq 0$ with ${\rho}_{ij}={\rho}^r_{i,rj}$ because we have ${\rho}^k_{l,ij}+{\rho}^k_{i,jl}+{\rho}^k_{j,li}=0$. With $\tau={\omega}^{ij}{\tau}^r_{ri,j}$ and $\rho={\omega}^{ij}{\rho}_{ij}$, we obtain $(n-2){\tau}^r_{ri,j}=(n-1){\rho}_{ij}+{\rho}_{ji}-(n/2(n-1)){\omega}_{ij}\rho$ and thus $n\rho=2(n-1)\tau$. The lower sequence splits with ${\varphi}_{ij}\rightarrow {\tau}_{ij}={\tau}^r_{ri,j}=(1/2){\varphi}_{ij}\rightarrow {\tau}_{ij}-{\tau}_{ji}={\varphi}_{ij}$ and ${\rho}_{ij}={\rho}_{ji}\Leftrightarrow {\varphi}_{ij}=0$ in $Z^2(g_1)\subset {\wedge}^2T^*\otimes g_1$. It follows from a chase that the kernel of the canonical projection $Riemann\rightarrow Weyl$ is defined by ${\rho}^k_{l,ij}={\tau}^k_{li,j}-{\tau}^k_{lj,i}$ with $({\rho}^k_{l,ij})\in Z^2(g_1)\subset Z^2({\hat{g}}_1)$ and $({\tau}^k_{li,j})\in T^*\otimes {\hat{g}}_2$. Accordingly $(n-2){\tau}_{ij}=n{\rho}_{ij}-(n/2(n-1)){\omega}_{ij}\rho$ provides the isomorphism $Ricci\simeq S_2T^*$ and we get $n{\rho}^k_{l,ij}={\delta}^k_i{\tau}_{lj}-{\delta}^k_j{\tau}_{li}+{\omega}_{lj}{\omega}^{ks}{\tau}_{si}-{\omega}_{li}{\omega}^{ks}{\tau}_{sj}$, that is: \\
\[ {\rho}^k_{l,ij}=  \frac{1}{(n-2)}({\delta}^k_i{\rho}_{lj}-{\delta}^k_j{\rho}_{li}+{\omega}_{lj}{\omega}^{ks}{\rho}_{si}-{\omega}_{li}{\omega}^{ks}{\rho}_{sj})-\frac{1}{(n-1)(n-2)}({\delta}^k_i{\omega}_{lj}-{\delta}^k_j{\omega}_{li})\rho       \]
We check that ${\rho}^r_{i,rj}={\rho}_{ij}$, obtaining therefore a splitting of the right vertical sequence in the last diagram that allows to define the Weyl tensor by difference. These purely algebraic results  only depend on $\omega$ independently of any conformal factor. \hspace{1cm}$ \Box $ \\

\noindent
{\bf EXAMPLE} 7.11: The free movement of a body in a constant static gravitational field $\vec{g}$ is described by $\frac{d\vec{x}}{dt}-\vec{v}=0, \frac{d\vec{v}}{dt}-\vec{g}=0, \frac{\partial\vec{g}}
{{\partial}x^i}-0=0$ where the "speed" is considered as a first order jet (Lorentz rotation) and the "gravity" as a second order jet (elation). Hence an {\it accelerometer} merely helps measuring the part of the Spencer operator dealing with second order jets ({\it equivalence principle}). As a byproduct, the difference ${\partial}_4f^k_4-f^k_{44}$ under the constraint ${\partial}_4f^k-f^k_4$ identifying the "speed" with a first order jet allows to provide a modern version of the {\it Gauss principle of least constraint} where the extremum is now obtained with respect to the second order jets and not with respect to the "acceleration" as usual ([1], p 470). The corresponding infinitesimal variational principle $\delta \int (\rho({\partial}_4{\xi}^4-{\xi}^4_4)+g^i({\partial}_i{\xi}^r_r-{\xi}^r_{ri})+g^{ij}({\partial}_i{\xi}^r_{rj}-0))dx=0$ provides the Poisson law of gravitation with $\rho=cst$ and $\vec{g}=(g^i)$ when ${g}^{ij}=\lambda{\omega}^{ij}\Rightarrow g_i=-{\partial}_i\lambda$. The last term of this gravitational action in vacuum is thus of the form $\lambda div(A)$, that is {\it exactly} the term responsible for the Lorentz constraint in Remark 7.6. \\

\noindent
{\bf 8   CONCLUSION}\\
   
\hspace*{5mm}In continuum mechanics, the classical approach is based on differential invariants and only involves derivatives of finite transformations. Accordingly, the corresponding variational calculus can only describe forces as it only involves translations. It has been the idea of E. and F. Cosserat to change drastically this point of view by considering a new differential geometric tool, now called Spencer sequence, and a corresponding variational calculus involving {\it both} translations and rotations in order to describe torsors, that is {\it both} forces and couples.\\
\hspace*{5mm}About at the same time, H. Weyl tried to describe electromagnetism and gravitation by using, {\it in a similar but complementary way}, the dilatation and elations of the conformal group of space-time. We have shown that the underlying Spencer sequence has additional terms, {\it not known today}, wich explain in a unique way all the above results and the resulting field-matter couplings.\\
\hspace*{5mm}In gauge theory, the structure of electromagnetism is coming from the unitary group $U(1)$, the unit circle in the complex plane, which is {\it not} acting on space-time, as the {\it only} possibility to obtain a pure 2-form from ${\wedge}^2T^*\otimes \cal {G}$ is to have $dim(G)=1$. However, we have explained the structure of electromagnetism from that of the conformal group of space-time, with a {\it shift by one step} in the interpretation of the Spencer sequence involved because the "{\it fields}" are now sections of $C_1\simeq T^*\otimes \cal{G}$ parametrized by $D_1$ and thus killed by $D_2$. \\
\hspace*{5mm}In general relativity, we have similarly proved that the standard way of introducing the Ricci tensor was based on a {\it double confusion} between the Janet and Spencer sequences described by {\it diagrams} (1) {\it and} (2). In particular we have explained why the intrinsic structure of this tensor {\it necessarily} depends on the difference existing between the Weyl group and the conformal group which is coming from second order jets, relating for the first time on equal footing electromagnetism and gravitation to the Spencer $\delta$-cohomology of various symbols.\\   
\hspace*{5mm}Accordingly, paraphrasing W. Shakespeare, we may say:\\
   
 \hspace*{25mm}         " TO  ACT  OR  NOT  TO  ACT,  THAT  IS  THE  QUESTION " \\
 
 \noindent
 and hope future will fast give an answer !.\\

 \noindent
 {    } \\
 
\noindent
{\bf BIBLIOGRAPHY}:\\

\noindent
[1] P. APPELL: Trait\'{e} de M\'{e}canique Rationnelle, Gauthier-Villars, Paris, 1909. Particularly t II concerned with analytical mechanics and t III with a Note by E. and F. Cosserat "Note sur la th\'{e}orie de l'action Euclidienne", 557-629.\\
\noindent
[2] V. ARNOLD: M\'{e}thodes math\'{e}matiques de la m\'{e}canique classique, Appendice 2 (G\'{e}od\'{e}siques des m\'{e}triques invariantes \`{a} gauche sur des groupes de Lie et hydrodynamique des fluides parfaits), MIR, moscow, 1974,1976. (For more details, see also: J.-F. POMMARET: Arnold's hydrodynamics revisited, AJSE-mathematics, 1, 1, 2009, 157-174).  \\
\noindent
[3] I. ASSEM: Alg\`{e}bres et Modules, Masson, Paris, 1997.\\
\noindent
[4] G. BIRKHOFF: Hydrodynamics, Princeton University Press, Princeton, 1954. French translation: Hydrodynamique, Dunod, Paris, 1955.\\
\noindent
[5] E. CARTAN: Sur une g\'{e}n\'{e}ralisation de la notion de courbure de Riemann et les espaces \`{a} torsion, C. R. Acad\'{e}mie des Sciences Paris, 174, 1922, 437-439, 593-595, 734-737, 857-860.\\
\noindent
[6] E. CARTAN: Sur les vari\'{e}t\'{e}s \`{a} connexion affine et la th\'{e}orie de la relativit\'{e} g\'{e}n\'{e}ralis\'{e}e, Ann. Ec. Norm. Sup., 40, 1923, 325-412; 41, 1924, 1-25; 42, 1925, 17-88.\\
\noindent
[7] O. CHWOLSON: Trait\'{e} de Physique (In particular III, 2, 537 + III, 3, 994 + V, 209), Hermann, Paris, 1914.\\
\noindent
[8] E. COSSERAT, F. COSSERAT: Th\'{e}orie des Corps D\'{e}formables, Hermann, Paris, 1909.\\
\noindent
[9] J. DRACH: Th\`{e}se de Doctorat: Essai sur une th\' {e}orie g\' {e}n\'{e}rale de l'int\'{e}gration et sur la classification des transcendantes, in Ann. Ec. Norm. Sup., 15, 1898, 243-384.\\
\noindent
[10] L.P. EISENHART: Riemannian Geometry, Princeton University Press, Princeton, 1926.\\
\noindent
[11] H. GOLDSCHMIDT: Sur la structure des \'{e}quations de Lie, J. Differential Geometry, 6, 1972, 357-373 and 7, 1972, 67-95.\\
\noindent
[12] H. GOLDSCHMIDT, D.C. SPENCER: On the nonlinear cohomology of Lie equations, I+II, Acta. Math., 136, 1973, 103-239.\\
\noindent
[13] M. JANET: Sur les syst\`{e}mes aux d\'{e}riv\'{e}es partielles, Journal de Math., 8, (3), 1920, 65-151. \\
\noindent
[14] E.R. KALMAN, Y.C. YO, K.S. NARENDA: Controllability of linear dynamical systems, Contrib. Diff. Equations, 1 (2), 1963, 189-213.\\
\noindent
[15] S. KOBAYASHI, K. NOMIZU: Foundations of Differential Geometry, Vol I, J. Wiley, New York, 1963, 1969.\\
\noindent
[16] G. KOENIG: Le\c{c}ons de Cin\'{e}matique (The Note "Sur la cin\'{e}matique d'un milieu continu" by E. Cosserat and F. Cosserat has rarely been quoted), Hermann, Paris, 1897, 391-417.\\
\noindent
[17] E.R. KOLCHIN: Differential Algebra and Algebraic Groups, Academic Press, New York, 1973.\\
\noindent
[18] A. KUMPERA, D.C. SPENCER: Lie Equations, Ann. Math. Studies 73, Princeton University Press, Princeton, 1972.\\
\noindent
[19] E. KUNZ: Introduction to Commutative Algebra and Algebraic Geometry, Birkh\"{a}user, 1985.\\
\noindent
[20] V. OUGAROV: Th\'{e}orie de la Relativit\'{e} Restreinte, MIR, Moscow, 1969; french translation, 1979.\\
\noindent
[21] W. PAULI: Theory of Relativity, Pergamon Press, London, 1958.\\
\noindent
[22] H. POINCARE: Sur une forme nouvelle des \'{e}quations de la m\'{e}canique, C. R. Acad\'{e}mie des Sciences Paris, 132 (7), 1901, 369-371.\\
\noindent
[23] J.-F. POMMARET: Systems of Partial Differential Equations and Lie Pseudogroups, Gordon and Breach, New York, 1978; Russian translation: MIR, Moscow, 1983.\\
\noindent
[24] J.-F. POMMARET: Differential Galois Theory, Gordon and Breach, New York, 1983.\\
\noindent
[25] J.-F. POMMARET: Lie Pseudogroups and Mechanics, Gordon and Breach, New York, 1988.\\
\noindent
[26] J.-F. POMMARET: Partial Differential Equations and Group Theory, Kluwer, Dordrecht, 1994.\\
\noindent
[27] J.-F. POMMARET: Partial Differential Control Theory, Kluwer, Dordrecht, 2001.\\
\noindent
[28] J.-F. POMMARET: Group interpretation of Coupling Phenomena, Acta Mechanica, 149, 2001, 23-39.\\
\noindent
[29] J.-F. POMMARET: Parametrization of Cosserat equations, Acta Mechanica, 215, 2010, 43-55.\\
\noindent
[30] J.-F. POMMARET: Macaulay inverse systems revisited, Journal of Symbolic Computation, 46, 2011, 1049-1069.\\
\noindent
[31] J.F. RITT: Differential Algebra, Dover, 1966.\\
\noindent
[32] J. J. ROTMAN: An Introduction to Homological Algebra, Academic Press, 1979.\\
\noindent
[33] D. C. SPENCER: Overdetermined Systems of Partial Differential Equations, Bull. Am. Math. Soc., 75, 1965, 1-114.\\
\noindent
[34] I. STEWART: Galois Theory, Chapman and Hall, 1973.\\
\noindent
[35] P.P. TEODORESCU: Dynamics of Linear Elastic Bodies, Editura Academiei, Bucuresti, Romania; Abacus Press, Tunbridge, Wells, 1975.\\
\noindent
[36] E. VESSIOT: Sur la th\'{e}orie des groupes infinis, Ann. Ec. Norm. Sup., 20, 1903, 411-451.\\
\noindent
[37] E. VESSIOT: Sur la th\'{e}orie de Galois et ses diverses g\'{e}n\'{e}ralisations, Ann. Ec. Norm. Sup., 21, 1904, 9-85.\\
\noindent
[38] C. YANG: Magnetic monopoles, fiber bundles and gauge fields, Ann. New York Acad. Sciences, 294, 1977, 86.\\
\noindent
[39] C.N. YANG, R.L. MILLS: Conservation of isotopic gauge invariance, Phys. Rev., 96, 1954, 191-195.\\
\noindent
[40] Z. ZOU, P. HUANG, Y. ZHANG, G. LI: Some researches on gauge theories of gravitation, Scientia Sinica, XXII, 6, 1979, 628-636.\\

\end{document}